\documentclass{amsart}
\usepackage{amssymb}
\usepackage[all]{xy}
\CompileMatrices
\newtheorem*{introthm}{Theorem}

\newtheorem{theorem}{Theorem}[section]
\newtheorem{lemma}[theorem]{Lemma}
\newtheorem{proposition}[theorem]{Proposition}
\newtheorem{corollary}[theorem]{Corollary}
\theoremstyle{definition}
\newtheorem{definition}[theorem]{Definition}
\newtheorem{example}[theorem]{Example}

\newtheorem{remark}[theorem]{Remark}

\numberwithin{equation}{theorem}
\def\div{{\rm div}}

\def\quot{/\!\!/}
\def\mal{\! \cdot \!}

\def\reg{{\rm reg}}
\def\rq#1{\widehat{#1}}

\def\b#1{\overline{#1}}
\def\bangle#1{\langle #1 \rangle}

\def\KK{{\mathbb K}}
\def\TT{{\mathbb T}}
\def\ZZ{{\mathbb Z}}

\def\QQ{{\mathbb Q}}
\def\PP{{\mathbb P}}

\def\WDiv{\operatorname{WDiv}}

\def\Cl{\operatorname{Cl}}

\def\Pic{\operatorname{Pic}}

\def\Supp{{\rm Supp}}
\def\Spec{{\rm Spec}}
\def\Proj{{\rm Proj}}

\def\Star{{\rm star}}

\def\cone{{\rm cone}}

\def\pr{{\rm pr}}

\newcounter{itemnumber}

\begin{document}
\title[GIT-equivalence beyond the ample cone]
{GIT-equivalence beyond the ample cone}
\author[F. Berchtold]{Florian Berchtold} 
\address{Fachbereich Mathematik und Statistik, Universit\"at Konstanz,
78457 Konstanz, Germany}
\email{florian.berchtold@uni-konstanz.de}
\author[J.~Hausen]{J\"urgen Hausen} 
\address{Mathematisches Institut, Universit\"at T\"ubingen,
Auf der Morgenstelle 10, 72076 T\"ubingen, Germany}
\email{juergen.hausen@uni-tuebingen.de}
\begin{abstract}
Given an algebraic torus action on a 
normal projective variety with finitely 
generated total coordinate ring, 
we study the GIT-equivalence for 
not necessarily ample linearized 
divisors, and we provide a 
combinatorial description
of the partially ordered set of 
GIT-equivalence classes.
As an application, we extend in
the $\QQ$-factorial case
a basic feature of the collection
of ample GIT-classes to 
the partially ordered collection of 
maximal subsets with a quasiprojective 
quotient: for any two members there 
is at most one minimal member 
comprising both of them.
Moreover, we demonstrate in an example, 
how our theory can be applied for a 
systematic treatment of 
``exotic projective orbit spaces'', 
i.e., projective geometric 
quotients that do not arise from any 
linearized ample divisor.
\end{abstract}

\maketitle

\section*{Introduction}

The approach to moduli spaces,
e.g., for curves of fixed genus,
presented by D.~Mumford in his 
Geometric Invariant Theory~\cite{Mu} 
relies on his
construction of quotients for 
actions of reductive groups $G$~on algebraic 
varieties~$X$. 
He introduces the notion of a 
$G$-linearized line bundle on $X$, and to any 
such bundle $L$ he associates a $G$-invariant open 
set $X^{ss}(L) \subset X$ of semistable points.
This set admits a so-called good quotient
$X^{ss}(L) \to  X^{ss}(L) \quot G$ with
a quasiprojective quotient space.

Mumford's construction, however, is in general not 
unique: his ``GIT-quotients'' turn out to depend 
essentially on the choice of the bundle and the 
linearization. 
Therefore, it is a natural desire to describe 
the collection of all possible GIT-quotients
for a given reductive group action. 
For ``ample GIT-quotients'', i.e.,
those arising from linearized ample line
bundles, this problem has meanwhile been studied
by several authors, see~\cite{BriPro}, 
\cite{DoHu}, \cite{Thadd}, and, finally, 
\cite{Re}.

A first basic step in the study of 
ample GIT-quotients is to show that 
there are only finitely many of them, 
see~\cite{DoHu}, \cite{Thadd}, 
\cite{BB}, and~\cite{AS}.
Then the subject becomes combinatorial. 
The situation is described
by sort of a fan subdividing the so-called 
(open) $G$-ample cone: the cones of this fan  
correspond to the ample GIT-quotients, 
and the face relations reflect in an order 
reversing manner the set theoretical 
inclusion of the respective sets of semistable 
points, see~\cite{Re}.

However, there are interesting examples of 
projective
GIT-quotients that do not arise from linearized 
ample bundles, see~\cite{BBSw2}.
Motivated by this observation, we study 
here the situation beyond the $G$-ample 
cone, and we propose a combinatorial 
framework for the description of the 
phenomena occuring there.
We restrict to the case of a torus action. 
On the one hand, concerning variation of 
GIT-quotients, 
this case is the most vivid one,
and, on the other, it allows an elementary 
treatment.

\goodbreak

The setup is the following:
$X$ is a normal projective variety
over an algebraically closed field 
$\KK$ of characteristic zero
such that $X$ has a free finitely 
generated divisor class group $\Cl(X)$, 
and a finitely generated total
coordinate ring 
(see Section~\ref{sec:semistabcrit})
$$
\mathcal{R}(X)
\;  := \; 
\bigoplus_{\Cl(X)} \Gamma(X,\mathcal{O}(D)).
$$
We consider the action $T \times X \to X$
of an algebraic torus $T = \Spec(\KK[M])$,
where $M$ is the lattice of characters. 
This comprises subtorus actions 
on projective toric varieties~\cite{KaStZe}, 
and, more generally, on projective 
spherical varieties.

In our setup, we can even do a little 
more than only describing 
the collection of the quotients arising
from the possible $T$-linearized line bundles over 
$X$: we allow
more generally quotients arising 
from $T$-linearized Weil divisors, 
see Section~\ref{sec:goodquot}.
Compared to Mumford's original approach,
this has the advantage that also for 
singular $X$, we obtain all good
quotients with a quasiprojective quotient
space, see~\cite{Ha2}.

We make use of the fact that $X$ is  
a good quotient of an invariant open 
subset $\rq{X}$ 
of $\b{X} := \Spec(\mathcal{R}(X))$
by the torus $H := \Spec(\KK([\Cl(X)])$
corresponding to the divisor class group.
The action of $T$ may be lifted to the 
multicone $\b{X}$ over $X$, and this lifting
corresponds to a refinement of the grading
$$ 
\mathcal{R}(X)
\; = \; 
\bigoplus_{(D,w) \in \Cl(X) \oplus M} \Gamma(X,\mathcal{R})_{(D,w)}.
$$
It turns out that the degrees 
$(D,w) \in \Cl(X) \oplus M$ 
of the refined grading
are in one-to-one
correspondence with the possible $T$-linearizations 
of the divisor classes $D \in \Cl(X)$, and that
the sets of semistable points only 
depend on the classes of the $T$-linearized divisors.
 
Let us indicate how the combinatorial description
runs. 
We associate to any $T$-linearized
class $(D,w)$ 
having a nonempty set of
semistable points
what we call its {\em GIT-bag\/}
$\mu(D,w)$.
This GIT-bag is a 
certain pointed convex polyhedral 
cone living in the rational vector space 
associated to $\Cl(X) \oplus M$, 
and it can be directly computed from a finite set 
of orbit data associated to the lifted action
of $H \times T$ on $\b{X}$.  
The set of GIT-bags is finite, and it comes
with a natural partial ordering ``$\le$''.
Our main result is the following, see
Theorem~\ref{GIT-bags}:

\begin{introthm}
Let $(D_i,w_i) \in \Cl(X) \oplus M$ represent 
two $T$-linearized Weil divisors on $X$,
and let $\mu(D_i,w_i)$ denote the 
associated GIT-bags. Then, for the 
associated sets of semistable points, 
we have 
$$
X^{ss}(D_1,w_1) \subset X^{ss}(D_2,w_2)
\iff
\mu(D_1,w_1) \ge \mu(D_2,w_2). $$
\end{introthm}

Inside the $T$-ample cone, the GIT-bags coincide
with the cones of the fan subdivision 
defined by the GIT-chambers of~\cite{DoHu}
and~\cite{Re}.
But outside the $T$-ample cone not much 
is left from the fan properties;
for example, overlappings are possible,
see Section~\ref{sec:examples}.
Nevertheless, the GIT-bags allow
to formulate answers to several 
questions.

For example, motivated by~\cite{BB}, 
we study  
{\em qp-maximal $T$-sets\/}. These are 
$T$-invariant open
subsets $U \subset X$ admitting a good
quotient $U \to U \quot T$ with a 
quasiprojective quotient space $U \quot T$
that do not occur as a saturated 
subset w.r. to the quotient map of 
a properly larger $U' \subset X$ with
the same properties.
Any qp-maximal set is a set of semistable
points of a linearized Weil divisor,
and in terms of GIT-bags, qp-maximality is 
characterized as follows, see 
Corollary~\ref{qp-maximal}:

\begin{introthm}
A GIT-bag describes a qp-maximal $T$-set
if and only if its relative interior 
is set theoretically minimal in the collection 
of the relative interiors of all GIT-bags.
\end{introthm}

Using this characterization, one 
can easily produce examples of
qp-maximal $T$-sets having a 
non-complete quotient
space, see Example~\ref{maxnotproj}.
It can as well be described in 
terms of GIT-bags, when a quotient
space is projective, see
Proposition~\ref{projectivequotients},
and there is a simple criterion
to figure out the ample GIT-quotients,
see Proposition~\ref{amplequotients}.
These two criteria are useful to
discuss an ``exotic orbit space'' 
presented in~\cite{BBSw2}, see 
Example~\ref{exoticquotients}.

For the case of a $\QQ$-factorial
variety $X$, the combinatorial
description of the collection 
of qp-maximal $T$-sets allows
to extend a basic statement 
from the ample theory: there,
one obtains as a consequence of the 
fan structure that any two 
sets of semistable points arising from
ample bundles admit at most one 
minimal such set comprising both
of them. 
We show the following,
see Corollary~\ref{minimum}.

\begin{introthm}
Suppose that $X$ is $\QQ$-factorial.
Then, for any two qp-maximal
$T$-sets $U_1, U_2 \subset X$,
the collection of qp-maximal $T$-sets
of $U \subset X$ with 
$(U_1 \cup U_2) \subset U$ is either
empty, or it contains a unique 
minimal element.  
\end{introthm}

The paper is organized as follows.
In Section~\ref{sec:goodquot},
we recall some basics on good quotients,
and the construction presented in~\cite{Ha2}.
Moreover, we introduce the group
of isomorphism classes of $G$-linearized
Weil divisors.
In Section~\ref{sec:affinecase}, we 
present a simple direct proof for the 
affine version of~\cite{Re}, 
which is needed later,
but also might be of independent
interest.
Section~\ref{sec:semistabcrit} is devoted to 
a combinatorial characterization of
semistability, and the main results
are presented in Section~\ref{sec:generalcase}.
Finally, in Section~\ref{sec:Qfactorialcase},
we investigate the case of a $\QQ$-factorial
$X$, and in Section~\ref{sec:examples},
we discuss some examples.
 
\tableofcontents

\section{Good quotients}
\label{sec:goodquot}

In this section, we recall 
the notion of a good quotient,
and we provide the basic facts
on this concept.
Moreover, we briefly recall 
from~\cite{Ha2} a generalization 
of Mumford's construction of
good quotients,
using Weil divisors instead of 
line bundles.
Finally, we introduce the group of
isomorphism classes of linearized
Weil divisors
and the GIT-equivalence.

\goodbreak

We work over an algebraically 
closed field $\KK$ of characteristic 
zero, the word variety refers to
a reduced scheme of finite type 
over $\KK$, and by a point we 
always mean a closed point.
When we speak of an action of an
algebraic group $G$ on a variety
$X$, we tacitly assume that this
action is given by a morphism
$G \times X \to X$; we then 
also speak of the $G$-variety $X$.

\begin{definition}
\label{goodquotdef}
Let a reductive linear algebraic group $G$ 
act on a variety $X$.
\begin{enumerate}
\item 
A {\em good quotient\/} for the $G$-action is
a $G$-invariant affine morphism $\pi \colon X \to Y$ 
such that the canonical map 
$\mathcal{O}_Y \to \pi_*(\mathcal{O}_X)^G$
is an isomorphism.
\item 
A good quotient $\pi \colon X \to Y$ is called 
{\em geometric\/}, if each fibre $\pi^{-1}(y)$,
where $y \in Y$, consists of a single $G$-orbit.
\end{enumerate}
\end{definition}

The definition of a good quotient was formulated 
by Seshadri~\cite{Se}, but, implicitly, the concept
occured already in Mumford's book~\cite{Mu}. 
Note that good quotients are obtained by glueing
classical Invariant Theory quotients $X \to Y$
of affine $G$-varieties~$X$, that means that
$Y := \Spec(\Gamma(X, \mathcal{O})^G)$ holds.
Here comes a list of basic properties, see for 
example~\cite{Se}:

\begin{proposition}
\label{goodquotprop}
Let $\pi \colon X \to Y $ be a good quotient for a
$G$-action on a variety~$X$.
\begin{enumerate}
\item
If $A,B \subset X$ are $G$-invariant closed 
subsets with $A \cap B = \emptyset$, then 
their images in $Y$ are as well closed and
satisfy $\pi(A) \cap \pi(B) = \emptyset$.
\item
For every $y \in Y$, there is a unique closed $G$-orbit 
$G \mal x_0 \subset \pi^{-1}(y)$, 
and this orbit lies in the closure of any other
orbit $G \mal x \subset \pi^{-1}(y)$.
\item 
If $\varphi \colon X \to Z$ is a $G$-invariant morphism,
then there is a unique morphism $\psi \colon Y \to Z$
with $\varphi = \psi \circ \pi$. 
\end{enumerate} 
\end{proposition}

The last property implies that a good quotient
for a $G$-variety is basically unique, provided
it exists. 
This justifies the notations $X \to X \quot G$ 
for a good quotient,
and $X \to X/G$ for a geometric quotient.
In general, a $G$-variety $X$ need not admit
a good quotient, but it may have many 
$G$-invariant open subsets $U \subset X$
with a good quotient $U \to U \quot G$.
For the study of such subsets, the 
following concept is crucial, 
compare~\cite{BB}.

\begin{definition}
\label{Gsatdef}
Let $X$ be a $G$-variety. 
A $G$-invariant open 
subset $U \subset X$ is called 
{\em $G$-saturated in $X$\/}
if we have $\b{G \mal x} \subset U$
for all $x \in U$,
where $\b{G \mal x}$ denotes the
orbit closure taken in $X$.
\end{definition}

Usually, one compares invariant open 
subsets $V \subset U$ 
of a $G$-variety $X$, that means that
one asks if $V$ is $G$-saturated in 
the $G$-variety
$U$.
If $G$ is reductive linear algebraic, and there is 
a good quotient $\pi \colon U \to U \quot G$,
then, by Property~\ref{goodquotprop}~(ii),
the set $V$ is $G$-saturated in $U$ if and only 
if $V = \pi^{-1}(\pi(V))$ holds.
In that case, $\pi(V)$ is open in $U \quot G$,
and the restriction 
$\pi \vert_V \colon V \to \pi(V)$
is a good quotient.

We now recall the construction
of good quotients given
in~\cite{Ha2}.
It extends Mumford's construction
by taking Weil divisors instead
of line bundles.
The advantage of this approach is 
that it provides also in the 
singular (normal) case
basically all good quotients
with a quasiprojective quotient. 
We should note that we present here 
a slightly modified version,
allowing also nontrivial linearizations
of the trivial divisor $D = 0$.
However, on the results and their
proofs, this has no impact.

\goodbreak

Let $X$ be a normal $G$-variety, 
where $G$ is a reductive linear
algebraic group. 
To any Weil divisor $D$ on $X$, 
we associate a sheaf of 
$\mathcal{O}_X$-algebras, 
and consider the corresponding
relative spectrum with its 
canonical morphism:
$$
\mathcal{A}
\; := \; 
\bigoplus_{n \in \ZZ_{\ge 0}}
\mathcal{O}_X(nD),
\qquad
X(D) \; := \;
\Spec_X(\mathcal{A}),
\qquad
q_D \colon X(D) \to X.
$$
The $\ZZ_{\ge 0}$-grading of the 
sheaf of algebras $\mathcal{A}$ 
defines an action of the multiplicative
group $\KK^* = \Spec(\KK[\ZZ])$ on 
$X(D)$, and the canonical morphism 
$q_D \colon X(D) \to X$ is a good 
quotient for this action.

Note that near singular points of
$X$, the scheme $X(D)$ need a priori
not be of finite type over $X$; however, 
we need not care to much about this 
difficulty, because in the situation 
we are interested in, it does not 
occur by assumption. 
Similarly to~\cite{Ha2}, 
we define:

\begin{definition}
A {\em $G$-linearization\/} of the divisor 
$D$ is a morphical $G$-action on $X(D)$ 
that commutes with the $\KK^*$-action on $X(D)$
and makes $q_D \colon X(D) \to X$ into a 
$G$-equivariant morphism.
\end{definition}

Any $G$-linearization of the divisor $D$ gives rise 
to a rational $G$-representation on the global 
sections respecting the $\ZZ_{\ge 0}$-grading, 
namely
$$
G \times \Gamma(X(D), \mathcal{O}) 
\to 
\Gamma(X(D),\mathcal{O}),
\qquad
(g \mal f)(x) := f(g^{-1} \mal x).
$$ 
Similarly to~\cite{Ha2}, we now introduce 
a notion of semistability. 
As usual, we denote for a section 
$f \in \Gamma(X,\mathcal{O}(D))$
of a Weil divisor its set of zeroes as
$$
Z(f)
\; := \;
\Supp(\div(f)+D). 
$$
\begin{definition}
\label{semistabdef}
Let $D$ be a $G$-linearized Weil divisor
on $X$. We call a point $x \in X$ 
{\em semistable\/} with respect to this
linearization, if there are 
an $n \in \ZZ_{>0}$ and a section
$f \in \Gamma(X,\mathcal{O}(nD))$ 
such that $X \setminus Z(f)$ is an
affine neighbourhood of $x$ and $f$ is 
invariant under the $G$-representation
on $\Gamma(X(D),\mathcal{O})$.
\end{definition}

We denote the set of semistable points
of a $G$-linearized Weil divisor $D$ 
by $X^{ss}(D)$, or by $X^{ss}(D,G)$
if we want to specify the group $G$.
From~\cite{Ha2} we infer the basic 
features of this construction:

\begin{proposition}
\label{qpquotients}
Let $G$ be a reductive linear algebraic group
and let $X$ be a normal $G$-variety.
\begin{enumerate}
\item If $D$ is a $G$-linearized Weil divisor on $X$,
then there exists a good quotient 
$X^{ss}(D) \to X^{ss}(D) \quot G$ with a 
quasiprojective quotient space. 
\item If $U \subset X$ is a $G$-invariant
open subset having a good quotient 
$U \to U \quot G$ with $U  \quot G$ 
quasiprojective, then $U$ is $G$-saturated
in some set $X^{ss}(D)$.
\end{enumerate}
\end{proposition}

In the literature, one often introduces 
a $G$-linearization of a line bundle
$L \to X$ over a $G$-variety more 
geometrically as a fibrewise linear
lifting of the $G$-action to the total
space $L$, see e.g.~\cite{KKV}.
Here comes the relation to our 
definition:

\begin{remark}
\label{dualrep}
If $D$ is Cartier and represents 
the class of a line bundle $L \to X$ in $\Pic(X)$,
then $X(D) \to X$ is the dual bundle of 
$L \to X$, and there is an isomorphism 
$$ 
\Gamma(X,L) \to \Gamma(X(D),\mathcal{O})_1,
\qquad
s \mapsto f_s,
\quad
\text{where } f_s(z) := \bangle{z,s(q_D(z))}.
$$
If $D$ is $G$-linearized, then the $G$-action
on $X(D)$ defines a dual, fibrewise 
linear action on the total space $L$ via
$$
\bangle{z, g \mal y}
\; := \; 
\bangle{g^{-1} \mal z, y}
\qquad
\text{ for }
g \in G, \; 
y \in L_x, \; 
z \in X(D)_{g \mal x}, \;
x \in X. 
$$
This action makes the projection
equivariant, and it induces the
so called ``dual representation''
of $G$ on the space  $\Gamma(X,L)$
of global sections: 
$$
g \mal s (x)
\; = \; 
g \mal (s(g^{-1} \mal x)).
$$ 
With respect to this  
representation, the isomorphism 
$\Gamma(X,L) \to \Gamma(X(D),\mathcal{O})_1$
mentioned before becomes an isomorphism
of $G$-modules.
\end{remark}

We conclude the section with a few words
about the passage to divisor classes.
For any $G$-variety $X$, there is the notion
of the group $\Pic_G(X)$ of isomorphism classes of
$G$-linearized line bundles over $X$.
Let us show how this concept can be extended 
to Weil divisors.

First of all, we have to prepare the definition 
of the $G$-linearized sum 
of two $G$-linearized Weil
divisors $D_1$ and $D_2$ on a normal 
$G$-variety $X$.
For this, we consider the
following sheaf of bigraded 
$\mathcal{O}_X$-algebras and its relative 
spectrum:
$$
\mathcal{B}
\; := \;
\bigoplus_{(n,m) \in \ZZ^2_{\ge 0}}
\mathcal{O}(nD_1+mD_2),
\qquad
X(D_1,D_2)
\; := \; 
\Spec_X(\mathcal{B}).
$$
Note that $X(D_1,D_2)$ comes with an
action of the torus 
$\TT^2 := \KK^* \times \KK^*$
defined by the bigrading of $\mathcal{B}$.
Moreover, we have canonical morphisms
$X(D_1,D_2) \to X(D_i)$ 
arising from the 
inclusions of sheaves
\begin{eqnarray*}
\bigoplus_{k \in \ZZ_{\ge 0}}
\mathcal{O}_X(kD_i)
& \to &
\bigoplus_{(n,m) \in \ZZ^2_{\ge 0}}
\mathcal{O}(nD_1+mD_2).
\end{eqnarray*}
These morphisms determine a morphism
$\varphi \colon X(D_1,D_2) \to X(D_1) \times_X X(D_2)$ 
to the fibre product,
which also comes with a canonical 
$\TT^2$-action and the diagonal 
$G$-action. 
Here are the basic properties 
of this setting.

\begin{lemma}
\label{sumlin}
For two $G$-linearized Weil divisors 
$D_1$, $D_2$ on $X$,
let $X(D_1,D_2)$ etc. be as above. 
Then there is a commutative diagram of 
$\TT^2$-equivariant morphisms:
$$
\xymatrix@!R0{
&&&
&& && && 
X(D_1) \times_X X(D_2) 
\ar[ddll]
\ar[drr]
\ar[ddd]
&& 
\\
&&&
X(D_1,D_2)
\ar@/^1pc/[urrrrrr]^{\varphi}
\ar[drrrr]
\ar[rrrrrrrr]|!{[rrrrr]}\hole|!{[rrrrrr]}\hole
\ar@/_1pc/[ddrrrrrr]_{q_{D_1,D_2}}
&& && && && 
 X(D_2)
\ar[ddll]
\\
&&&
&& && 
X(D_1)  
\ar[drr]
&& && && && 
\\
&&&
&& && && X &&
}
$$
The diagonal $G$-action on the fibre product lifts 
uniquely to $X(D_1,D_2)$, and then descends 
to $X(D_1+D_2)$ within a further canonical 
commutative diagram: 
$$
\xymatrix{
X(D_1,D_2) \ar[rr] \ar[dr]
& & 
X(D_1+D_2) \ar[dl]
\\
& X &
}
$$
Moreover, the induced $G$-action on $X(D_1+D_2)$
is a $G$-linearization of the Weil divisor $D_1+D_2$. 
\end{lemma}

\begin{proof}
The morphism $\varphi \colon X(D_1,D_2) \to
X(D_1) \times_X X(D_2)$ is given by the universal
property of the fibre product.
It is an isomorphism over the set $X_{\reg} \subset X$ 
of smooth points, because there it comes from 
the canonical isomorphism of the corresponding 
sheaves of $\mathcal{O}_X$-algebras:
\begin{eqnarray*}
\left(
\bigoplus_{m \in \ZZ_{\ge 0}}
\mathcal{O}(mD_1)
\right)
\otimes_{\mathcal{O}_X}
\left(
\bigoplus_{n \in \ZZ_{\ge 0}}
\mathcal{O}(nD_2)
\right)
& \to &
\bigoplus_{(m,n) \in \ZZ^2_{\ge 0}}
\mathcal{O}(mD_1+nD_2).
\end{eqnarray*}
Moreover, as this is a bigraded homomorphism, 
we can conclude that $\varphi$ is 
$\TT^2$-equivariant.
The fact that $\varphi$ is an isomorphism
over the $G$-invariant set $X_{\reg} \subset X$,
allows us to shift the diagonal $G$-action 
on the fibre product over $X_{\reg}$ to a 
morphical action
$$ 
\alpha \colon
G \times q_{D_1,D_2}^{-1}(X_{\reg})
\; \to \;  
q_{D_1,D_2}^{-1}(X_{\reg}).
$$
Our task is to extend this action to the whole 
$X(D_1,D_2)$. This is done via extending the 
corresponding comorphisms. 
Let $\beta \colon G \times X \to X$
denote the action on $X$.
Then, by $G$-equivariance, we obtain
a commutative diagram  
$$ 
\xymatrix{
{\Gamma(U \cap X_{\reg},\mathcal{B})}
\ar@{=}[d]
\ar[rr]^{\alpha^*}
& & 
{\Gamma(\beta^{-1}(U \cap X_{\reg}), \mathcal{O}_G \otimes \mathcal{B})}
\ar@{=}[d]
\\
{\Gamma(U,\mathcal{B})}
\ar[rr]
& & 
{\Gamma(\beta^{-1}(U), \mathcal{O}_G \otimes \mathcal{B})}
}
$$
for any affine open subset $U \subset X$.
As one easily verifies, the lower rows of 
these diagrams are the 
comorphisms of a $G$-action on $X(D_1,D_2)$.
By construction, this extension has the desired
properties.
So, the first part of the lemma is proved.

To see the second part,
consider the antidiagonal $\KK^*$-action on 
$X(D_1,D_2)$
defined by the homomorphism of tori 
$\KK^* \to \TT^2$ sending 
$t$ to $(t,t^{-1})$.
This action admits a good quotient, namely
the morphism $X(D_1,D_2) \to X(D_1 + D_2)$
arising from the canonical injection of
sheaves
\begin{eqnarray*}
\bigoplus_{k \in \ZZ_{\ge 0}}
\mathcal{O}(kD_1+kD_2)
& \to &
\bigoplus_{(m,n) \in \ZZ^2_{\ge 0}}
\mathcal{O}(mD_1+nD_2).
\end{eqnarray*}
Since the antidiagonal $\KK^*$-action and the 
$G$-action on $X(D_1,D_2)$ commute, 
the $G$-action descends to an
action on the quotient
space $X(D_1+D_2)$. 
By construction, it commutes
with the $\KK^*$-action on 
 $X(D_1+D_2)$, 
and the canonical morphism
$X(D_1+D_2) \to X$ becomes $G$-equivariant.
\end{proof}

\begin{definition}
Let $D_1$ and $D_2$ be two $G$-linearized 
Weil divisors on a normal $G$-variety $X$.
\begin{enumerate}
\item
The $G$-linarization of the sum
$D_1 + D_2$ is the unique $G$-action on 
$X(D_1+D_2)$ provided by Lemma~\ref{sumlin}.
\item 
We say that $D_1$ and $D_2$ are isomorphic,
if there is a $(\KK^* \times G)$-equivariant
isomorphism $X(D_1) \to X(D_2)$ over~$X$.
\end{enumerate}
\end{definition}

Note that for the case of a pair
of linearized Cartier
divisors, our definition of the 
linearized sum corresponds to  
the usual tensor product of 
linearized line bundles,
and the notion of isomorphism
is the usual one.
%The crucial properties of these 
%settings are the following:

\begin{proposition}
The set of isomorphism classes of
$G$-linearized Weil divisors form
a group $\Cl_G(X)$. Moreover,
\begin{enumerate}
\item forgetting about the linearizations
gives rise to a well defined homomorphism
$\Cl_G(X) \to \Cl(X)$,
\item
for any $G$-linearized Weil divisor
$D$ on $X$, the set $X^{ss}(D,G)$
only depends on its class in
$\Cl_G(X)$.
\end{enumerate}
\end{proposition}

Observe that the kernel of the forgetting
homomorphism $\Cl_G(X) \to \Cl(X)$
consists precisely of the linearizations 
of the trivial bundle.
Finally, by the above proposition, 
we may generalize 
the usual concept of GIT-equivalence to
the setting of Weil divisors.

\begin{definition}
We say that two $G$-linearized divisor
classes in $\Cl_G(X)$ are
{\em GIT-equivalent\/}
if they define the same set of semistable
points.
\end{definition}

\section{The affine case}
\label{sec:affinecase}

In this section, we study the collection 
of sets of semistable points 
arising from the possible linearizations 
of the trivial bundle over an affine
variety with a torus action. 
We provide a simple direct proof
for the fact that this collection 
is in order 
reversing bijection to a (quasi-)fan 
subdividing 
the weight cone of the action.

This result may be viewed as an affine 
version of~\cite{Re}.
It is well known for linear
torus actions on $\KK^n$; 
in this case, the describing fan is
a so-called 
Gelfand-Kapranov-Zelevinsky
decomposition, see~\cite{OdPa}.

The precise setup is the following. 
$\KK$ is an algebraically closed field,
$R$ is a finitely generated integral
$\KK$-algebra, graded by a lattice
$M \cong \ZZ^d$:
$$
R
\; = \; 
\bigoplus_{w \in M} R_w .
$$
This grading corresponds
to an action of the 
algebraic torus
$T := \Spec(\KK[M])$ on the 
affine variety $X := \Spec(R)$.

We denote by $M_{\QQ} := M \otimes_{\ZZ} \QQ$
the rational vector space associated to $M$.
Recall that the {\em weight cone\/} of 
the $T$-variety $X$ is the convex cone 
in $M_{\QQ}$ generated by all $w \in M$
admitting a nontrivial homogeneous $f \in R_w$:
$$
\Omega_T(X) 
\; := \;
\cone(w \in M; \; R_w \ne 0)
\; \subset \;
M_{\QQ}.
$$
Since the algebra $R$ is generated by finitely 
many homogeneous elements, 
the weight cone $\Omega_T(X)$ is finitely generated 
as well, and thus it is a polyhedral cone.
Note that $\Omega_T(X)$ is pointed, 
i.e., contains no line,  if  
$R_0 = \KK$ and $R^* = \KK^*$ hold. 

\begin{definition}
For a point $x \in X$, its {\em orbit monoid\/}
is the semigroup consisting of
all weights that admit a homogeneous 
function, which is invertible near $x$:
$$ 
S_{T}(x) 
\; := \; 
\{w \in M; \; \exists \ f \in R_w, \; f(x) \ne 0\}.
$$ 
The {\em orbit cone\/} of $x \in X$ 
is the convex (polyhedral) cone 
$\omega_T(x) \subset M_{\QQ}$ 
generated by the orbit monoid $S_T(x)$.
\end{definition}

We collect some basic properties
of orbit cones.
A first observation is that
the orbit cones are not affected by 
passing to the normalization.

\begin{lemma}
\label{normalization}
Let $\nu \colon X' \to X$ be the $T$-equivariant
normalization. Then, for every $x' \in X'$,
we have $\omega_T(x') = \omega_T(\nu(x'))$.  
\end{lemma}

\begin{proof}
The inclusion $\omega_T(\nu(x')) \subset \omega_T(x')$
is clear by equivariance.
The reverse inclusion follows from considering
equations of integral dependence for 
the homogeneous elements 
$f \in \mathcal{O}(X')$
with $f(x') \ne 0$. 
\end{proof}

We shall use the orbit cones to
describe properties of orbit closures.
The basic statement in this regard 
is the following one.

\begin{proposition}
\label{orbitcones1}
For a point $x \in X$, let $T_x \subset T$
be its isotropy group, and let
$M_T(x) \subset M$ be the sublattice
generated by the orbit monoid $S_T(x)$.
\begin{enumerate}
\item The algebraic torus 
$T/T_x$ acts with a dense free
orbit on the orbit closure 
$\b{T \mal x} \subset X$.
\item 
$\b{T \mal x}$ 
has the affine toric variety 
$\Spec(\KK[\omega_T(x) \cap M_T(x)])$
as its $(T/T_x)$-equivariant normalization.
\end{enumerate}
\end{proposition}

\begin{proof}
The first assertion is obvious,
and the second one follows immediately
from Lemma~\ref{normalization} and 
the fact  
that the algebra of global functions 
of $\b{T \mal x}$ is the semigroup
algebra $\KK[S_T(x)]$ of the weight
monoid. 
\end{proof}

For two polyhedral cones $\omega_1$ and 
$\omega_2$ in a common vector space
we write $\omega_1 \preceq \omega_2$ 
if $\omega_1$ is a face of $\omega_2$.
Lemma~\ref{normalization} and 
Proposition~\ref{orbitcones1} have the 
following consequence.

\begin{corollary}
\label{clos2face}
Let $x \in X$. Then the $T$-orbits in 
$\b{T \mal x}$ correspond to the faces
of $\omega_T(x)$ via
$T \mal y \mapsto \omega_T(y)$.
\end{corollary}

The following simple observation will 
replace in our setup the deeper finiteness 
result on GIT-quotients given in~\cite{DoHu}
and~\cite{Thadd}.

\begin{proposition}
\label{orbitcones2}
The collection of orbit cones 
$\{\omega_T(x); \; x \in X\}$ 
is finite.
\end{proposition}

\begin{proof}
Embed $X$ equivariantly into some $\KK^n$, 
on which $T$ acts diagonally. 
Then the $T$-orbit cone of a point $x \in X$
equals its $T$-orbit cone w.r. to $\KK^n$.
The $T$-orbit cones w.r. to $\KK^n$
are constant along the orbits of the
standard action of $\TT^n := (\KK^*)^n$,
because this action commutes with 
that of $T$.
Since $\TT^n$ has only finitely 
many orbits in $\KK^n$, 
the assertion follows.
\end{proof}

We now enter the study of
the collection of sets of semistable
points arising from the possible 
$T$-linearizations of the trivial bundle.
First of all, let us recall that
these linearizations correspond to 
the characters of $T$.

\begin{lemma}
\label{lem:lintriv}
Consider a $T$-linearization of the trivial
bundle over $X$. Then there is a unique 
$w \in M$ such that the dual $T$-action 
on $X \times \KK$ is of the form
\begin{eqnarray}
\label{eq:lintriv}
t \mal (x,z)
& := &
(t \mal x, \chi^w(t)z). 
\end{eqnarray}
\end{lemma}

\begin{proof}
The dual action is a fibrewise linear $T$-action
on $X \times \KK$ making 
$X \times \KK \to X$ equivariant. 
Consequently, there is a morphism 
$c \colon T \times X \to \KK^*$ 
such that we have
\begin{eqnarray*}
t \mal (x,z)
& := &
(t \mal x, c(t,x) z).
\end{eqnarray*}
Clearly, we always have $c(1,x) = 1$.
Thus, for fixed $x$, the map $t \mapsto c(t,x)$ 
is a homomorphism.
Hence, by rigidity of tori, $c$ does not depend on
$x$.
\end{proof}

In the sequel we shall denote by 
$X^{ss}(w) \subset X$ 
the set of semistable points defined
by the linearization~\ref{eq:lintriv}. 
It can be 
explicitly described in terms of 
homogeneous functions and also 
in terms of orbit cones.

\begin{lemma}
\label{orbitcones2semistab}
The set $X^{ss}(w) \subset X$
of semistable points of the 
linearization~\ref{eq:lintriv}
is given by
\begin{eqnarray*}
X^{ss}(w)
& = &
\bigcup_{f \in R_{nw}, \; n \in \ZZ_{>0}} 
X_f 
\\
& & 
\\
& = & 
\{x \in X; \; w \in \omega_T(x)\}.
\end{eqnarray*}
In particular, the set of semistable points 
$X^{ss}(w)$ is nonempty 
if and only if $w \in \Omega_T(X) \cap M$ holds.
\end{lemma}

\begin{proof}
As indicated in Remark~\ref{dualrep}
the invariant sections for the 
linearization~\ref{eq:lintriv} are 
precisely the functions $f \in R_{nw}$ with
$n \in \ZZ_{\ge 0}$.
This gives the first equality.
The second one is a direct consequence
of the definition of an orbit cone,
and the last statement is obvious. 
\end{proof}

As outlined in Section~\ref{sec:goodquot}, 
the set $X^{ss}(w)$ is 
$T$-invariant, and it admits a good quotient 
$X^{ss}(w) \to Y(w)$ by the action of $T$. 
In fact, the quotient space
$Y(w) = X^{ss}(w) \quot T$ is the homogeneous 
spectrum of a Veronese subalgebra:
$$ 
Y(w) 
\; = \; 
\Proj(R(w)),
\qquad
\text{where} \quad
R(w) 
\; = \; \bigoplus_{n \in \ZZ_{\ge 0}} R_{nw}
\; \subset \; 
R.  
$$
In particular, every quotient space $Y(w)$ 
is projective over $Y(0) = \Spec(R_0)$. 
Furthermore, if we have $X^{ss}(w_1) \subset X^{ss}(w_2)$, 
then Proposition~\ref{goodquotprop}~(iii)
yields a commutative diagram
\begin{equation*}
\label{eq:gitsyst}
\xymatrix{
X^{ss}(w_1) \ar[rr]^{\subset}  \ar[d]_{\quot T}  & & X^{ss}(w_2)
\ar[d]^{\quot T} \\ 
Y(w_1) \ar[rr]^{{\varphi^{w_1}_{w_2}}} \ar[rd] & & Y(w_2) \ar[ld] \\
& Y(0)
}
\end{equation*}
Note that the induced map 
$\varphi^{w_1}_{w_2} \colon Y(w_1) \to Y(w_2)$ 
of the quotient spaces is dominant and projective.
Moreover, we have 
$\varphi^{w_1}_{w_3} 
=  
\varphi^{w_2}_{w_3} \circ \varphi^{w_1}_{w_2}$
whenever composition is possible.

The collection of all nonempty sets $X^{ss}(w)$
together with their good quotients $X^{ss}(w) \to Y(w)$ 
and the above diagrams 
is called the {\em GIT-system\/} associated 
to the trivial bundle on the $T$-variety $X$. 
Let us turn to the combinatorial description 
of this GIT-system.
We introduce a collection of convex, polyhedral 
cones.

\begin{definition}
For a weight $w \in \Omega_T(X) \cap M$, 
the associated {\em GIT-cone\/} 
is the (nonempty) intersection of all orbit 
cones containing $w$:
$$ 
\sigma_{T}(w) 
\; := \; 
\bigcap_{w \in \omega_T(x)} \omega_T(x).
$$ 
Moreover, 
the collection of all the possible
GIT-cones defined by the action of $T$
on $X$ is denoted as
$$
\Sigma_T(X) 
\; := \; 
\{\sigma_T(w); \; w \in \Omega_T(X) \cap M\}.
$$
\end{definition}

Note that, for us, GIT-cones are closed cones,
and thus they are not chambers in the sense 
of~\cite{Re}.
A first important observation is that the GIT-cones
are in order reversing one-to-one correspondence 
with the possible 
sets of semistable points arising from the various
linearizations of the trivial bundle.

\begin{proposition}
\label{interior}
Let $w_1,w_2 \in \Omega_T(X) \cap M$.
Then we have  
\begin{enumerate}
\item
$ 
\displaystyle
X^{ss} (w_1) \subset X^{ss} (w_2)
\iff 
\sigma_T(w_1) \supset \sigma_T(w_2),
$

\smallskip

\item
$
\displaystyle 
X^{ss} (w_1) = X^{ss} (w_2)
\iff 
\sigma_T(w_1) = \sigma_T(w_2).
$
\end{enumerate}  
\end{proposition}

\begin{proof}
This is an immediate consequence of the definition
of the GIT-cones and the characterization of
semistability in terms of orbit cones given in
Lemma~\ref{orbitcones2semistab}
\end{proof}

This proposition allows us to speak about the set of 
semistable points corresponding to a GIT-cone 
$\sigma \in \Sigma_T(X)$: we set 
$$ 
X^{ss}(\sigma) 
\; := \; 
X^{ss}(w),
\text{ where } 
\sigma = \sigma_T(w).
$$

\begin{lemma}
\label{chamber2semistabset}
The set of semistable points associated to a 
GIT-cone $\sigma \in \Sigma_T(X)$ is given 
by
$$ 
X^{ss}(\sigma)
\; = \; 
\{x \in X; \; \sigma \subset \omega_T(x)\}.
$$
\end{lemma}

We now come to the main result of this section.
Together with Proposition~\ref{interior}, it 
describes the structure of the collection of sets
of semistable points associated to the 
linearizations of the trivial bundle
as a partially ordered set.

A {\em quasifan\/} 
is a finite collection $\Sigma$ 
of  (not necessarily pointed) convex,
polyhedral cones in a common vector space
such that for $\sigma \in \Sigma$ also 
all faces of $\sigma$ belong to $\Sigma$, 
and for any two 
$\sigma, \sigma' \in \Sigma$
the intersection 
$\sigma \cap \sigma'$ 
is a face of both, 
$\sigma$ and $\sigma'$.
A quasifan is called a {\em fan\/}
if it consists of pointed cones.
The {\em support\/} of a quasifan is 
the union of its cones.

\begin{theorem}
\label{affineGIT}
The collection of all GIT-cones
$\Sigma_T(X)$
is a quasifan in the vector space 
$M_{\QQ}$ having the weight
cone $\Omega_T(X)$ as its support. 
\end{theorem}

In the proof of this result, we
need a further basic property of  
the GIT-cones, also needed 
later.
For a convex polyhedral cone $\sigma$, we denote 
its relative interior by $\sigma^{\circ}$,
that means that $\sigma^{\circ}$ is obtained
by removing all proper faces from $\sigma$. 

\begin{lemma}
\label{intchamber}
Let $w \in \Omega_T(X) \cap M$.
Then the associated GIT-cone
$\sigma := \sigma_T(w) \in \Sigma_T(X)$
satisfies
$$
\begin{array}{ccccccc} 
&  
& \sigma 
&  = 
&  \displaystyle 
\bigcap_{w \in \omega_T(x)^{\circ}} \omega_T(x)
&  = 
& \displaystyle 
\bigcap_{\sigma^{\circ} \subset \omega_T(x)^{\circ}} \omega_T(x),
\\
\\
w 
& \in 
& \sigma^{\circ} 
&  =  
& \displaystyle 
\bigcap_{w \in \omega_T(x)^{\circ}} \omega_T(x)^{\circ}
& =  
&\displaystyle 
\bigcap_{\sigma^{\circ} \subset \omega_T(x)^{\circ}} \omega_T(x)^{\circ}.
\end{array}
$$
\end{lemma}

\begin{proof}
For any orbit cone $\omega_T(x)$ with
$w \in \omega_T(x)$, there is a unique 
minimal face $\omega \preceq \omega_T(x)$
with $w \in \omega$. This face satisifies
$w \in \omega^{\circ}$. According to 
Corollary~\ref{clos2face}, the face 
$\omega \preceq \omega_T(x)$ is again 
an orbit cone. 
This gives the first formula.
The second one follows from an elementary
observation: if the intersection 
of the relative interiors of a finite number 
convex polyhedral cones is nonempty, then it 
equals the relative interior of the intersection 
of the cones. 
\end{proof}

\begin{proof}[Proof of Theorem~\ref{affineGIT}]
First of all note that, by finiteness 
of the number of orbit cones, as shown 
in Proposition~\ref{orbitcones2}, 
the collection of all GIT-cones is 
finite.

The further proof is split   
into verifications of several claims.
For the sake of handy notation, we
set for the moment 
$\Omega := \Omega_T(X)$ and 
$\Sigma := \Sigma_T(X)$.
Moreover, we omit the 
subscript ``$T$''
when denoting orbit cones 
and GIT-cones,
and we write $X(\sigma)$
instead of $X^{ss}(\sigma)$.

\medskip

\noindent
{\em Claim~1. } 
Let $\sigma_1,\sigma_2 \in \Sigma$ 
with $\sigma_1 \subset \sigma_2$.
Then, for every $x_1 \in X(\sigma_1)$
with 
$\sigma_1^{\circ} 
\subset 
\omega(x_1)^{\circ}$,
there exists an $x_2 \in X(\sigma_2)$ 
with 
$\omega(x_1) \preceq \omega(x_2)$.

\medskip

Let us verify the claim.
By Proposition~\ref{interior}, we have
$X(\sigma_2) \subset X(\sigma_1)$.
Consequently, the GIT-system 
provides a commutative diagram
with a dominant, proper, hence 
surjective morphism 
$\varphi \colon Y(\sigma_2) \to Y(\sigma_1)$
of the quotient spaces:
$$
\xymatrix{
X(\sigma_2) 
\ar[r]^{\subset}
\ar[d]_{p_2}
&
X(\sigma_1) 
\ar[d]^{p_1}
\\
Y(\sigma_2) 
\ar[r]_{\varphi} 
&
Y(\sigma_1) 
}
$$
If a point $x_1 \in X(\sigma_1)$ satisfies
$\sigma_1^{\circ} \subset \omega_T(x_1)^{\circ}$,
then, by Lemma~\ref{chamber2semistabset} 
and Corollary~\ref{clos2face}, 
its $T$-orbit is closed in $X(\sigma_1)$.
Proposition~\ref{goodquotprop}~(ii) thus tells us
that $x_1 \in \b{T \mal x_2}$ holds for any 
point $x_2$ belonging to the (nonempty) intersection
$X(\sigma_2) \cap p_1^{-1}(p_1(x_1))$.
Using once more Corollary~\ref{clos2face}
gives Claim~1.

\medskip

\noindent
{\em Claim~2. } 
Let $\sigma_1,\sigma_2 \in \Sigma$. 
Then $\sigma_1 \subset \sigma_2$
implies $\sigma_1 \preceq \sigma_2$.

\medskip

For the verification, let 
$\tau_2 \preceq \sigma_2$
be the (unique) face with 
$\sigma_1^{\circ} \subset \tau_2^{\circ}$,
and let  
$\omega_{1,1}, \ldots, \omega_{1,r}$
be the orbit cones with 
$\sigma_1^{\circ} \subset \omega_{1,i}^{\circ}$.
Then we obtain, using Lemma~\ref{intchamber} for the
second observation, 
$$
\tau_2^{\circ} \cap \omega_{1,i}^{\circ} \; \ne \; \emptyset,
\qquad
\sigma_1
\; = \; 
\omega_{1,1} \cap \ldots \cap \omega_{1,r}.
$$
By Claim~1, we have 
$\omega_{1,i} \preceq \omega_{2,i}$
with orbit cones 
$\omega_{2,i}$ satisfying
$\sigma_2 \subset \omega_{2,i}$,
and hence $\tau_2 \subset \omega_{2,i}$.
The first of the displayed formulas 
implies $\tau_2 \subset \omega_{1,i}$,
and the second one thus gives $\tau_2 = \sigma_1$.
So, Claim~2 is verified.

\medskip

\noindent
{\em Claim~3. } 
Let $\sigma \in \Sigma$. 
Then every face 
$\sigma_0 \preceq \sigma$ 
belongs to $\Sigma$.

\medskip

To see this, consider any 
$w \in \sigma_0^{\circ}$.
Lemma~\ref{intchamber} yields 
$w \in \sigma(w)^{\circ}$.
By the definition of GIT-cones, 
we have $\sigma(w) \subset \sigma$. 
Claim~2 gives even $\sigma(w) \preceq \sigma$.
Thus, we have two faces,
$\sigma_0$ and $\sigma(w)$ of $\sigma$
having a common point $w$ in their relative
interiors.
This means $\sigma_0 = \sigma(w)$,
and Claim~3 is verified.

\medskip

\noindent
{\em Claim~4. } 
Let $\sigma_1, \sigma_2 \in \Sigma$. 
Then $\sigma_1 \cap \sigma_2$ 
is a face of both, 
$\sigma_1$ and $\sigma_2$.

\medskip

Let $\tau_i \preceq \sigma_i$ be
the minimal face containing 
$\sigma_1 \cap \sigma_2$. 
Choose $w$ in the relative interior 
of $\sigma_1 \cap \sigma_2$,
and consider the GIT-cone
$\sigma(w)$.
By Lemma~\ref{intchamber} and the 
definition of GIT-cones, we see
$$ 
w \; \in \; 
\sigma(w)^{\circ} \cap \tau_i^{\circ},
\qquad
\sigma(w) 
\; \subset \;
\sigma_1 \cap \sigma_2 
\; \subset \; 
\tau_i.
$$
By Claim~2, the second relation implies
in particular
$\sigma(w) \preceq \sigma_i$.
Hence, we can conclude 
$\sigma(w) = \tau_i$, and
hence $\sigma_1 \cap \sigma_2$
is a face of both $\sigma_i$. 
Thus, Claim~4 is verified,
and the properties of a quasifan are
established for $\Sigma_T(X)$.
\end{proof}

\section{A semistability criterion}
\label{sec:semistabcrit}

We present a combinatorial description 
of the set of semistable points 
associated to a linearized Weil divisor.
By $X$ we denote a normal projective variety 
with finitely generated 
free divisor class group $\Cl(X)$,
and we consider the action 
$T \times X \to X$
of an algebraic torus
$T = \Spec(\KK[M])$ on 
the variety $X$.

The {\em total coordinate ring\/} $\mathcal{R}(X)$
of the variety $X$ is defined as follows:
choose a subgroup $K \subset \WDiv(X)$
of the group of Weil divisors such that
the canonical map $K \to \Cl(X)$ is an
isomorphism, and set
$$ 
\mathcal{R}(X)
\; := \;
\Gamma(X,\mathcal{R}), 
\qquad
\text{where} 
\quad
\mathcal{R}
\; := \;
\bigoplus_{D \in K} 
\mathcal{O}(D). 
$$
This ring depends only up to isomorphism
on the choices made in its definition.
An important property of the total coordinate
ring $\mathcal{R}(X)$ is that it is a factorial ring, 
see~\cite{BeHa1} and~\cite{ElKuWa}.

Throughout this section,
we assume that $\mathcal{R}(X)$ is 
finitely generated as a $\KK$-algebra.
We consider the following geometric
objects associated to the 
$K$-graded sheaf $\mathcal{R}$ of
$\mathcal{O}_X$-algebras:
$$
H \; := \; \Spec(\KK[K]),
\qquad
\b{X} := \Spec(\mathcal{R}(X)),
\qquad
\rq{X}
\; := \; 
\Spec_X(\mathcal{R}).
$$
So, $\rq{X}$ refers
to the relative spectrum of $\mathcal{R}$.
Recall that there is a canonical morphism
$q_X \colon \rq{X} \to X$.
We list some basic properties
of this setting,
which will be frequently used
in the subsequent constructions 
and proofs, compare also~\cite{BeHa3}.

\goodbreak

\begin{proposition}
\label{lift}
Let $H$, $\b{X}$, $\rq{X}$ and 
$q_X \colon \rq{X} \to X$ 
be as before. Then the following statements
hold.
\begin{enumerate}
\item
The $K$-grading of $\mathcal{R}$ defines 
an action of the torus $H$ on $\rq{X}$,
and $q_X \colon \rq{X} \to X$ is a 
good quotient for this action.
\item
The $K$-grading of $\mathcal{R}(X)$
defines an action of  the torus $H$ on $\b{X}$,
and the canonical map $\rq{X} \to \b{X}$ 
is an equivariant open embedding.
\item
For $D \in K$ and  
$f \in \Gamma(X,\mathcal{O}(D))$ 
with $X \setminus Z(f)$ affine,
the inverse image
$q_X^{-1}(X \setminus Z(f))$
equals $\b{X}_f$.
\item
For the set $X_{\reg} \subset X$ of 
nonsingular points, the complement
$\b{X} \setminus q_X^{-1}(X_{\reg})$
is of codimension at least two in
$\b{X}$.
\item
Suppose that $H \mal x \subset \rq{X}$ 
is closed.
Then $f \in \Gamma(\b{X},\mathcal{O})_D$
satisfies $f(x) = 0$ if and only if 
one has $q_X(x) \in Z(f)$ for 
$f \in \Gamma(X,\mathcal{O}(D))$.  
\item
There exists a $T$-action 
on $\b{X}$, commuting with the 
$H$-action on $\b{X}$
such that $\rq{X} \subset \b{X}$ 
is $T$-invariant and 
$q_X \colon \rq{X} \to X$ is 
$T$-equivariant.
\end{enumerate}
\end{proposition}

\begin{proof}
We begin with a basic observation.
Let $D \in K$ and 
$f \in \Gamma(X,\mathcal{O}(D))$ 
be such that $X \setminus Z(f)$ is 
affine.
Then there are the following 
identities of global
functions:
$$
\Gamma(\b{X}_f, \mathcal{O})
\; = \; 
\mathcal{R}(X)_f 
\; = \; 
\Gamma(X \setminus Z(f),\mathcal{R})
\; = \; 
\Gamma(q_X^{-1}(X \setminus Z(f)), \mathcal{O}).
$$
Since we have $K=\Cl(X)$, the
variety $X$ is covered by such
affine sets $X \setminus Z(f)$. 
Thus, we see
in particular that $\mathcal{R}$ 
is locally of finite
type and $\rq{X}$ is a variety.

The first assertion is then obvious.
In the second one, only the claim that
$\rq{X} \to \b{X}$ is an open embedding
needs some explanation.
By the above identities, each affine
subset $q_X^{-1}(X \setminus Z(f))$
is mapped isomorphically onto $\b{X}_f$.
It follows that $\rq{X} \to \b{X}$
is an open embedding.
Moreover, Assertion~(iii) 
drops out as well.

The fourth assertion is 
due to a further identity of global
functions: it follows from the fact 
that $\Gamma(X_{\reg},\mathcal{R})$ 
equals  $\Gamma(X,\mathcal{R})$.

To verify  Assertion~(v), 
suppose first that $q_X(x) \not\in Z(f)$
holds for a section $f \in \Gamma(X,\mathcal{O}(D))$. 
Then $f$ restricts to an invertible section 
of $\mathcal{R}$ over a suitable neighbourhood
$U \subset X$ of $q_X(x)$.
Consequently, $f$ is invertible as a function 
on $q_X^{-1}(U)$, which implies $f(x) \ne 0$.

Conversely, let $f(x) \ne 0$ 
for $f \in \Gamma(\b{X},\mathcal{O})_D$.
Consider 
the orbit $H \mal x$, and the zero set 
$B := N(f,\rq{X})$.
By Proposition~\ref{goodquotprop}~(i),
the image $q_X(B) \subset X$ is closed and
does not contain $q_X(x)$.
Hence, for a suitable neighbourhood
$U \subset X$ of $q_X(x)$, we see  
that $f$ is invertible as a function
on $q_X^{-1}(U)$, and hence it is so 
as a section of $\mathcal{R}$ over $U$.
This implies $q_X(x) \not\in Z(f)$. 

So, we are left with verifying the last
statement. 
By~\cite{Ha2}, there is a $T$-linearization 
of the group $K \subset \WDiv(X)$ 
over $X_{\reg} \subset X$.
In other words, we may lift the $T$-action
to $q_X^{-1}(X_{\reg})$.
By the first assertion, the complement
$\b{X} \setminus q_X^{-1}(X_{\reg})$
is of codimension at least two in $\b{X}$.
Hence the lifted $T$-action extends to
$\b{X}$.
\end{proof}

For the remainder of this section, 
we fix a lifting of the $T$-action 
to $\b{X}$ as in Proposition~\ref{lift}~(vi). 
In terms of multigraded rings, this means
that we have a refinement of the $K$-grading:
\begin{eqnarray*}
\mathcal{R}(X)
& = & 
\bigoplus_{(D,w) \in K \oplus M} \Gamma(X,\mathcal{R})_{(D,w)}.
\end{eqnarray*}

We need a pullback construction
for linearized Weil divisors.
For a Weil divisor $D$ on $X$, 
consider its restriction $D_{\reg}$
to $X_{\reg}$, and let
$\b{D}$ denote the Weil divisor on
$\b{X}$ obtained by closing the 
support of $q_X^*D_{\reg}$.
Now, suppose that $D$ is $T$-linearized.
We observe that then $\b{D}$ inherits 
in a canonical 
way an $(H \times T)$-linearization.
In fact, consider the cartesian 
square
$$
\xymatrix{
{q_X^{-1}(X_{\reg})(q_X^*D_{\reg})}
\ar[d]\ar[r]
& 
{X_{\reg}}(D_{\reg})
\ar[d]^{q_D}
\\
{q_X^{-1}(X_{\reg})}
\ar[r]_{q_X}
&
X_{\reg}
}
$$
Viewing the upper left
space as a fibre product
$q_X^{-1}(X_{\reg}) \times_{X_\reg}  
X_{\reg}(D_{\reg})$, 
one defines an $(H \times T)$-action 
on it by letting $H$ act on the first 
factor and letting 
$T$ act diagonally. 
Since $\b{X}$ is locally factorial,
$\b{X}(\b{D}) \to \b{X}$ is a bundle, 
and thus, by Proposition~\ref{lift}~(iv),
the $(H \times T)$-action 
extends to the desired 
linearization of $\b{D}$.

\begin{lemma}
\label{saturated2}
The set $q_X^{-1}(X^{ss}(D,T))$
is $(H \times T)$-saturated in 
${\b{X}^{ss}(\b{D}, H \times T)}$,
and there is a commutative 
diagram, where the horizontal arrows
are open embeddings:
$$
\xymatrix{
{\rq{X}} 
\ar[d]_{\quot H}
&
q_X^{-1}(X^{ss}(D,T)) 
\ar[l]
\ar[r] \ar[d]^{\quot H}
&
{\b{X}^{ss}(\b{D}, H \times T)}
\ar[r] \ar[dd]^{\quot H \times T}
&
{\b{X}} 
\\
X
&
X^{ss}(D,T) 
\ar[l] \ar[d]^{\quot T}
&
&
\\
&
X^{ss}(D,T) \quot T 
\ar[r]
&
{\b{X}^{ss}(\b{D}, H \times T) \quot H \times T}
&
}
$$
\end{lemma}

\begin{proof}
Every $T$-invariant section
$f \in \Gamma(X,\mathcal{O}(nD))$
defines via pullback an 
$(H \times T)$-invariant section
of $\Gamma(\b{X}, \mathcal{O}(n\b{D}))$. 
Thus, Proposition~\ref{lift}~(iii)
and the definition of semistability
give the desired statement. 
\end{proof}

\begin{lemma}
\label{Tlindescr}
As a $(H \times T)$-linearized
divisor, $\b{D}$ is isomorphic to the trivial
bundle with an $(H \times T)$-linearization,
and there is a unique $w \in M$ such that 
the corresponding dual action is given as 
\begin{eqnarray*}
(h,t) \mal (x,z) 
& = &
((h,t) \mal x, \chi^{(D,w)}(h,t)z). 
\end{eqnarray*}
Moreover, the assignment $D \mapsto (D,w)$ 
induces an isomorphism 
$\Cl_T(X) \to K \oplus M$ from the group
of $T$-linearized divisor classes on 
$X$ to the character lattice of
the torus $H \times T$.
\end{lemma}

\begin{proof}
Consider the set $X_{\reg} \subset X$ of 
nonsingular points, and the restriction 
$D_{\reg}$ of $D$ to $X_{\reg}$. 
Then, over the sets $U_i \subset X_{\reg}$ 
of a suitably fine open cover, the sheaf 
$\mathcal{O}(D_{\reg})$ is 
generated by invertible elements 
$f_i \in \Gamma(U_i,\mathcal{O}(D))$.

The line bundle $\pi \colon L \to X_{\reg}$ 
with the transistion functions 
$\xi_{ij} := f_j/f_i$ is the dual 
bundle of $X_{\reg}(D_{\reg}) \to X_{\reg}$;
it comes with the dual $T$-action
and with canonical trivializations
$$
\pi^{-1}(U_i) \to U_i \times \KK,
\qquad 
v \mapsto (\pi(v),z_i(v)).
$$

The pullback line bundle 
$q_X^*L = q_X^{-1}(X_{\reg}) \times_{X_{\reg}} L$
is dual to the line bundle
$q_X^{-1}(X_{\reg})(q_X^*D_{\reg}) \to q_X^{-1}(X_{\reg})$
arising from the restriction of $\b{D}$.
The dual $(H \times T)$-action on 
$q_X^* L$ equals the pullback linearization,
and is of the form
$$ 
(h,t) \mal (x,v)
= (t \mal h \mal x, t \mal v).
$$

We claim that $q_X^*L$ is $H$-equivariantly
isomorphic to the trivial bundle, $H$-linearized 
by the character $\chi^{D}$;
this follows from the fact that the 
functions $f_i \in \Gamma(U_i,\mathcal{O}(D))$
define a global trivialization for $q_X^*L$, 
namely 
$$
q_X^{-1}(X_{\reg}) \times_{X_{\reg}} L \to q_X^{-1}(X_{\reg}) \times \KK,
\quad
(x,v) \mapsto (x,f_i(x)z_i(v)),
\text{ for } 
x \in q_X^{-1}(U_i).
$$

Using Proposition~\ref{lift}~(iv),
we can extend this to a global trivialization
of the dual bundle of $\b{D}$.
This proves the first part of the 
assertion.

For the second part, note first
that $H$ acts freely on $q_X^{-1}(X_{\reg})$,
because locally on $X_{\reg}$ all divisors $D \in K$
are principal, and hence, locally any point 
of $q_X^{-1}(X_{\reg})$ has $K$ as its weight 
monoid.
Moreover there is a commutative
diagram 
$$
\xymatrix{
0 \ar[r]
& 
M \ar[r] 
&
{\Pic_{H \times T}(q_X^{-1}(X_{\reg}))} \ar[r] 
&
{\Pic_{H}(q_X^{-1}(X_{\reg}))} \ar[r] 
&
0 
\\
0 \ar[r]
& 
M \ar[r] \ar@{=}[u]
&
{\Pic_T(X_{\reg})} \ar[r] \ar[u]^{q_X^*}
&
{\Pic(X_{\reg})} \ar[r] \ar[u]_{q_X^*}
&
0
} 
$$
having exact rows.
Since $H$ acts freely on $q_X^{-1}(X_{\reg})$,
we infer from~\cite[Proposition~4.2]{KKV} 
that the right hand side pullback 
is an isomorphism.
Consequently, also the pullback in the middle 
of the above diagram must be an isomorphism.

The assertion thus follows from the fact 
that
we have canonical isomorphisms
$\Cl_T(X) \cong \Pic_T(X_{\reg})$ 
and 
$\Cl_{H \times T}(\b{X}) 
\cong 
\Pic_{H \times T}(q_X^{-1}(X_{\reg}))$,
where the latter relies on
Proposition~\ref{lift}~(iv).
\end{proof}

Via the isomorphism $[D] \mapsto (D,w)$
just established, we shall identify 
from now on the 
$T$-linearized divisor classes on 
$X$ with the elements on $K \oplus M$.
We denote the corresponding
sets of semistable points by
$$
X^{ss}(D,w)
\; := \; 
X^{ss}(D,T),
\qquad
\text{where }
\Cl_T(X) \ni [D]
\; \mapsto \;
(D,w) \in K \oplus M.
$$

We are now ready to begin with the 
combinatorial characterization of 
semistability.
It involves two fans; 
namely the collections of GIT-cones
$\Sigma_{H \times T}(\b{X})$ and 
$\Sigma_{H}(\b{X})$ for the 
actions of $H \times T$ and $H$ on
$\b{X}$.
These collections are actually fans, 
because, by~\cite[Proposition~4.3]{BeHa3},
the weight cones 
$\Omega_{H \times T}(\b{X})$ and 
$\Omega_H(\b{X})$ are pointed.

Let $\kappa_X \in \Sigma_H(\b{X})$ be the 
GIT-cone
corresponding to $\rq{X} \subset \b{X}$,
which, by projectivity of $X$, is a set of 
$H$-semistable points. 
Moreover, let $\Pi \colon K \oplus M \to K$
denote the projection.
We consider the following collections of 
orbit cones:
\begin{eqnarray*}
C_T(X)
& := & 
\{\omega_{H \times T}(x); \; 
x \in \b{X}, \; 
\kappa_X^{\circ} \subset \Pi(\omega_{H \times T}(x))^{\circ}\},
\\[.5em]
C_T(\sigma) 
& := &
\{\omega_{H \times T}(x); \; 
x \in \b{X}, \;  
\sigma^{\circ} \subset \omega_{H \times T}(x)^{\circ}\},
\end{eqnarray*}
where $\sigma \in \Sigma_{H \times T}(\b{X})$ 
may be any GIT-cone.
The geometric meaning of these collections is that
they describe the collection of closed orbits in 
the respective sets of semistable points:

\begin{lemma}
\label{geominterp}
Let $\sigma \in \Sigma_{H \times T}(\b{X})$,
let $x \in \b{X}$, and consider the
orbit cone $\omega_{H \times T}(x)$.
\begin{enumerate}
\item The orbit $H \mal x$ is a closed subset of $\rq{X}$
      if and only if 
      $\omega_{H \times T}(x) \in C_T(X)$ holds. 
\item The orbit $(H \times T) \mal x$ is a closed subset of 
      $\b{X}^{ss}(\sigma)$
      if and only if 
      $\omega_{H \times T}(x) \in C_T(\sigma)$
      holds. 
\end{enumerate}
\end{lemma}

\begin{proof}
First note that for any orbit cone 
$\omega_{H \times T}(x)$, the image
$\Pi(\omega_{H \times T}(x))$
equals the orbit cone 
$\omega_H(x)$.
Thus, the collections $C_T(X)$ and $C_T(\sigma)$
describe the orbits of 
$H$ in $\rq{X} = \b{X}^{ss}(\kappa_X)$
and $H \times T$ in $\b{X}^{ss}(\sigma)$ 
having minimal orbit cones.
The assertions hence follow from 
Corollary~\ref{clos2face}.
\end{proof}

Our characterization of semistability is formulated
in terms of the above collections of orbit cones:

\begin{theorem}
\label{semsitablecrit}
Let $(D,w) \in \Omega_{H \times T}(\b{X}) \cap (K \oplus M)$ 
represent a $T$-linearized Weil divisor $D$ on~$X$, 
and consider the GIT-cone
$\sigma := \sigma_{H \times T}(D,w)$
in $\Sigma_{H \times T}(\b{X})$.
Then we have
$$ 
q_X^{-1}(X^{ss}(D,w))
\; = \; 
\{x \in \b{X}; \; 
\omega \preceq \omega_{H \times T}(x) 
\text{ for some } \omega \in C_T(X) \cap C_T(\sigma)\}.
$$
\end{theorem}

\begin{proof}
As before, let $\b{D}$ denote the divisor
on $\b{X}$, obtained by closing the support of
the pullback divisor $q_X^*D_{\reg}$, 
and consider the inherited
$(H \times T)$-linearization of~$\b{D}$.
Then $\b{X}^{ss}(\b{D}, H \times T)$ is 
precisely the set of semistable points 
corresponding to the GIT-cone 
$\sigma \in \Sigma_{H \times T}(\b{X})$.

To verify the ``$\subset$''-part,
consider first a closed orbit 
$(H \times T) \mal x$ 
in $q_X^{-1}(X^{ss}(D,T))$.
By Lemma~\ref{saturated2}, 
this orbit is closed in  
$\b{X}^{ss}(\b{D}, H \times T)$.
Thus, Lemma~\ref{geominterp} yields 
$\sigma^{\circ} \subset \omega_{H \times T}(x)^{\circ}$.
Moreover, the orbit $H \mal x$
is closed in $(H \times T) \mal x$,
hence in $q_X^{-1}(X^{ss}(D,T))$,
and thus, since $q_X^{-1}(X^{ss}(D,T))$ 
is $H$-saturated in $\rq{X}$,
it is even closed in $\rq{X}$.
Lemma~\ref{geominterp} yields  
$
\kappa_X^{\circ} \subset \omega_{H}(x)^{\circ} 
$.
So, 
$\omega_{H \times T}(x)$
lies in 
$
C_T(X) \cap
C_T(\sigma)
$.

Now, given an arbitrary point
$x \in q_X^{-1}(X^{ss}(D,T))$,
we may consider any point $y$ in the 
$(H \times T)$-orbit closure of
$x$ having a closed 
$(H \times T)$-orbit in 
$q_X^{-1}(X^{ss}(D,T))$.
According to Corollary~\ref{clos2face},
the orbit cone 
$\omega := \omega_{H \times T}(y)$ 
is a face of 
$\omega_{H \times T}(x)$
and,
by the preceding consideration,
$\omega$ belongs to $C_T(X) \cap C_T(\sigma)$.

We turn to the inclusion ``$\supset$''.
First consider a point 
$x \in \b{X}$ such that
$\omega_{H \times T}(x)$ belongs 
to $C_T(X) \cap C_T(\sigma)$.
Then we have 
$$ 
x 
\; \in \; \rq{X} \cap \b{X}^{ss}(\b{D}, H \times T).
$$
Moreover, by Lemma~\ref{geominterp},
the orbit $H \mal x$ is closed in $\rq{X}$, 
and the orbit $(H \times T) \mal x$
is closed in $\b{X}^{ss}(\b{D}, H \times T)$.

By a repeated shrinking procedure,
we shall now construct a neighbourhood
of $q_X(x) \in X$ as required in 
Definition~\ref{semistabdef}.
First, note that the definition of semistability 
for the linearized divisor $\b{D}$ provides an 
$f \in \Gamma(\b{X},\mathcal{O})$,
homogeneous of weight $(nD,nw)$
with some $n \in \ZZ_{> 0}$, such that
we have
$$
(H \times T) \mal x
\; \subset \; 
\b{X}_f 
\; \subset \; 
\b{X}^{ss}(\b{D}, H \times T).
$$

Consider the complement 
$B_1 := \b{X}_f \setminus \rq{X}$.
This is an 
$(H \times T)$-invariant closed 
subset of $\b{X}_f$ disjoint from
$(H \times T) \mal x$.
By Proposition~\ref{goodquotprop}~(i),
the good quotient
$$\b{X}_f \to \b{X}_f \quot (H \times T)$$ 
separates $x$ and $B_1$.
Thus, we can choose an 
$(H \times T)$-invariant function
$f_0 \in \Gamma(\b{X}_f,\mathcal{O})$
satisfying $f_0 \vert_{B_1} = 0$ and 
having no zeroes in $(H \times T) \mal x$.

For a suitable $k \in \ZZ_{>0}$, the product
$g := f_0f^k$ is a $T$-invariant 
element of $\Gamma(X, \mathcal{O}(knD))$.
Since $H \mal x$ is closed in $\rq{X}$,
Proposition~\ref{lift}~(v) yields
$$
x 
\; \in \;  
q_X^{-1}(X \setminus Z(g)) 
 \; \subset \; 
\b{X}_g
 \; \subset \;
\rq{X}.
$$

Now, consider the intersection
$B_2 := \b{X}_g \cap q_X^{-1}(Z(g))$.  
This is an $(H \times T)$-invariant 
closed subset of $\b{X}_g$ 
disjoint from 
$(H \times T) \mal x$.  
Similarly as before, we can choose an 
$(H \times T)$-invariant function
$g_0 \in \Gamma(\b{X}_g,\mathcal{O})$
satisfying $g_0 \vert_{B_2} = 0$ and 
having no zeroes in $(H \times T) \mal x$.

Once more, 
for a suitable $l \in \ZZ_{>0}$, the product
$h := g_0g^l$ is a $T$-invariant 
element of $\Gamma(X, \mathcal{O}(lknD))$.
This time we have 
$$
x 
\; \in \;  
q_X^{-1}(X \setminus Z(h)) 
 \; \subset \; 
\b{X}_h
\; \subset \;
q_X^{-1}(X \setminus Z(g)) 
\; \subset \;
\rq{X}.
$$

We claim that even $q_X^{-1}(X \setminus Z(h)) = \b{X}_h$
holds.
Assume, to the contrary, that there exists 
a point $y \in \b{X}_h$ with $q_X(y) \in Z(h)$.
Note that the orbit closure in $\rq{X}$ 
satisfies
$$ 
\b{H \mal y}
\; \subset \;
q_X^{-1}(q_X(y))
\; \subset \; 
q_X^{-1}(X \setminus Z(g)). 
$$
Consider any $y_0 \in \b{H \mal y}$ such that 
$H \mal y_0$ is closed in $\rq{X}$.
By the above observation,
$q_X(y_0) \not\in Z(g)$ holds. 
Proposition~\ref{lift}~(v) thus yields
$g(y_0) \ne 0$.
By assumption, we have 
$$
q_X(y_0)
\; = \; 
q_X(y) 
\; \in \; 
Z(h).
$$
Applying again Proposition~\ref{lift}~(v)
gives $h(y_0) = 0$, and thus $g_0(y_0) = 0$.
Since $g_0$ is an $H$-invariant function, this
means $g_0(y) = 0$. Thus, we obtain $h(y) = 0$,
which is in contradiction to $y \in \b{X}_h$.

Having seen that $q_X^{-1}(X \setminus Z(h)) = \b{X}_h$
holds, we easily obtain the rest: the element
$h \in \Gamma(X, \mathcal{O}(lknD))$
is $T$-invariant and defines an 
affine neighbourhood
$$ 
X \setminus Z(h)
\; = \; 
\b{X}_h \quot H
$$ 
of the point $q_X(x) \in X$ as required 
in Definition~\ref{semistabdef}.
This shows that the point $x$ 
belongs to $q_X^{-1}(X^{ss}(D,T))$.

If $x \in \b{X}$ is an arbitrary point 
belonging to the right hand side
set of the equation in the assertion, 
then Corollary~\ref{clos2face} tells us
that the face 
$\omega \preceq \omega_{H \times T}(x)$
with $\omega \in C_T(X) \cap C_T(\sigma)$ 
is the orbit cone of some 
point $y$ belonging to the 
$(H \times T)$-orbit closure 
of $x$. 
By the preceding consideration,
we have $y \in  q_X^{-1}(X^{ss}(D,T))$.
This implies $x \in  q_X^{-1}(X^{ss}(D,T))$.
\end{proof}

Even in the case of a trivial torus 
action, Theorem~\ref{semsitablecrit}
is of some interest:
it then provides a description of the 
cone of ample divisors of the variety 
$X$, compare also~\cite[Theorem~7.3]{BeHa3}.

\begin{corollary}
\label{amplecone}
The cone of ample divisor classes on
$X$ is the relative interior 
$\kappa_X^{\circ} \subset 
K_{\QQ} = \Cl_{\QQ}(X)$
of the GIT-cone $\kappa_X \in \Sigma_H(\b{X})$.
\end{corollary}

\begin{proof}
Consider the action of the trivial torus
$T = \{e_T\}$ on $X$.
Then, for any $x \in \b{X}$, we have 
$\omega_{H \times T}(x) = \omega_{H}(x)$.
Moreover, the fans $\Sigma_{H \times T}(\b{X})$
and $\Sigma_H(\b{X})$ coincide.
Any divisor $D \in K$ is $T$-linearized,
and $D \in K$ is ample if and only if
$X^{ss}(D) = X$ holds.
The latter is equivalent to 
$q_X^{-1}(X^{ss}(D))=\rq{X}$,
and, by Theorem~\ref{semsitablecrit},
this holds if and only
if we have $D \in \kappa_X^{\circ}$.
\end{proof}

\begin{remark}
The case of a trivial $T$-action already shows
that $q_X^{-1}(X^{ss}(D,T))$ is in general 
properly smaller than 
$\rq{X} \cap \b{X}^{ss}(\b{D}, H \times T)$.
Let $D$ be effective but not big. 
Then $X^{ss}(D,T)$ is empty, but 
$\rq{X} \cap \b{X}^{ss}(\b{D}, H \times T)$
is nonempty.
As an explicit example, one may take 
$X = \PP_1 \times \PP_1$ and 
$D = \PP_1 \times \{0\}$.  
\end{remark}

\section{The general case}
\label{sec:generalcase}

In this section, we present the main 
results of the paper.
As in the preceding section,
$X$ is a normal projective variety
with finitely generated total
coordinate ring $\mathcal{R}(X)$,
and the algebraic torus 
$T = \Spec(\KK[M])$ acts on $X$.
We give a combinatorial
description of the collection of sets
of semistable points associated
to the $T$-linearized Weil divisors
on $X$.

Recall from the preceding section
that $X$ is a good quotient
of an open subset $\rq{X}$ of 
the affine variety 
$\b{X} = \Spec(\mathcal{R}(X))$
by the torus $H = \Spec(\KK[K])$
corresponding to the grading 
lattice $K \cong \Cl(X)$ of
$\mathcal{R}(X)$.
As before, we fix a lifting of 
the $T$-action to $\b{X}$;
this corresponds to the choice
of a refined grading
$$ 
\mathcal{R}(X)
\; = \; 
\bigoplus_{(D,w) \in K \oplus M} \Gamma(X,\mathcal{R})_{(D,w)}.
$$

As observed in Lemma~\ref{Tlindescr},
the degrees $(D,w) \in K \oplus M$ 
describe the possible $T$-linearizations 
of the divisors $D \in K$.
Again, we denote by 
$\kappa_X \in \Sigma_H(\b{X})$ 
the GIT-cone corresponding to
$\rq{X} \subset \b{X}$.
Moreover, $\Pi \colon K \oplus M \to K$ 
denotes the projection, and we use the 
collection of cones
\begin{eqnarray*}
C_T(X)
& := & 
\{\omega_{H \times T}(x); \; 
x \in \b{X}, \; 
\kappa_X^{\circ} \subset \Pi(\omega_{H \times T}(x))^{\circ}\}.
\end{eqnarray*}

We first have to figure out the linearized 
divisor classes with a nonempty set of semistable
points. 
For that purpose, consider the set 
$$
C_T^{\sharp}(X)
\; := \; 
\bigcup_{\omega \in C_T(X)} \omega^{\circ}
\; \subset \; 
K_{\QQ} \oplus M_{\QQ}.
$$

\begin{lemma}
The set 
$C_T^{\sharp}(X)$ is a convex cone
in $K_{\QQ} \oplus M_{\QQ} = \Cl_{T}(X)_{\QQ}$, 
and for a vector $(D,w) \in K \oplus M$,
the set $X^{ss}(D,w)$ is nonempty
if and only if we have
$(D,w) \in C_T^{\sharp}(X)$.
\end{lemma}

\begin{proof}
The fact that $X^{ss}(D,w)$ is nonempty
if and only if
$(D,w) \in C_T^{\sharp}(X)$ holds, 
follows directly from  
Lemma~\ref{saturated2} and
Theorem~\ref{semsitablecrit}.
Moreover,
by multiplying suitable invariant sections,
one sees that for any two linearized divisor 
classes $(D_i,w_i)$ with nonempty sets of 
semistable points $X^{ss}(D_i,w_i)$, 
also $(D_1+D_2,w_1+w_2)$ 
admits semistable points.
This gives convexity of $C_T^{\sharp}(X)$.
\end{proof}

\begin{definition}
Let the pair $(D,w) \in C_T^{\sharp}(X) \cap (K \oplus M)$
represent a $T$-linearized divisor on $X$.  
Then we define its associated {\em GIT-bag\/} 
to be the convex polyhedral cone
$$ 
\mu(D,w)
\; = \; 
\bigcap_{\omega \in C_T(X); \; (D,w) \in \omega^{\circ}}
\omega.
$$ 
The collection of all these GIT-bags is denoted by
$\Lambda_T(X)$. 
For $\mu_1, \mu_2 \in \Lambda_T(X)$, we write 
$\mu_1 \le \mu_2$ if for any $\omega_2 \in C_T(X)$ 
with $\mu_2^{\circ} \subset \omega_2^{\circ}$ 
there is a face $\omega_1 \preceq  \omega_2$
with $\omega_1 \in C_T(X)$ and  
$\mu_1^{\circ} \subset \omega_1^{\circ}$.
\end{definition}

Note that every GIT-bag is a union of GIT-cones 
of the GIT-fan $\Sigma_{H \times T}(\b{X})$ 
in $K \oplus M$ corresponding 
to the $(H \times T)$-action on $\b{X}$.
Moreover, the relation ``$\le$'' 
clearly is a partial ordering on  
$\Lambda_T(X)$.
We shall now see that the partially ordered 
set of GIT-bags describes precisely 
the GIT-equivalence:

\begin{theorem}
\label{GIT-bags}
Let $(D_i,w_i) \in  C_T^{\sharp}(X) \cap (K \oplus M)$ 
represent two $T$-linearized Weil divisors on $X$.
Then we have 
$$
X^{ss}(D_1,w_1) \subset X^{ss}(D_2,w_2)
\iff
\mu(D_1,w_1) \ge \mu(D_2,w_2). $$
\end{theorem}

\begin{proof}
We shall make repeated use of the 
combinatorial characterization of
semistability given in 
Theorem~\ref{semsitablecrit}.
For this, let 
$\sigma_1, \sigma_2 \in \Sigma_{H \times T}(\b{X})$ 
denote the GIT-cones
associated to  
$(D_1,w_1)$ and $(D_2,w_2)$, respectively.
Moreover, set
$$ 
W_i 
\; := \; 
q_X^{-1}(X^{ss}(D_i,w_i))
\; \subset \; 
\b{X}.
$$

First suppose that 
$X^{ss}(D_1,w_1) \subset X^{ss}(D_2,w_2)$
holds.
Consider an orbit cone 
$\omega_1 = \omega_{H \times T}(x) \in C_T(X)$ 
with $\mu(D_1,w_1) \subset \omega_1^{\circ}$.
Then $(D_1,w_1) \in \omega_1^{\circ}$
and, by Lemma~\ref{intchamber}, 
$\sigma_1^{\circ} \subset \omega_1^{\circ}$
hold.
Theorem~\ref{semsitablecrit} yields 
$x \in W_1$. 
Thus, by assumption, 
$x \in W_2$ holds.
Again by Theorem~\ref{semsitablecrit},
there is an
$\omega_2 \preceq \omega_1$ 
with $\omega_2 \in C_T(X)$ and 
$\sigma_2^{\circ} \subset \omega_2^{\circ}$.
Thus, we see $(D_2,w_2) \in \omega_2^{\circ}$, and
hence, $\mu(D_2,w_2)^{\circ} \subset \omega_2^{\circ}$.
This eventually implies 
$\mu(D_1,w_1) \ge \mu(D_2,w_2)$. 

Conversely, suppose that  
$\mu(D_1,w_1) \ge \mu(D_2,w_2)$
holds. 
Consider $x \in W_1$
with $(H \times T)\mal x$ closed 
in $W_1$.
By Theorem~\ref{semsitablecrit}
and Corollary~\ref{clos2face}, 
the orbit cone 
$\omega_1 := \omega_{H \times T}(x)$ 
belongs to $C_T(X)$ and satisfies 
$\sigma_1^{\circ} \subset \omega_1^{\circ}$.
The latter implies 
$(D_1,w_1) \in \omega_1^{\circ}$,
and hence   
$\mu(D_1,w_1)^{\circ} \subset \omega_1^{\circ}$.
By assumption, there is an
$\omega_2 \preceq \omega_1$ with
$\omega_2 \in C_T(X)$ and 
$\mu(D_2,w_2)^{\circ} \subset \omega_2^{\circ}$.
The latter implies 
$\sigma_2^{\circ}\subset \omega_2^{\circ}$. 
Thus, Theorem~\ref{semsitablecrit}
yields $x \in W_2$,
and we can conclude $W_1 \subset W_2$.
\end{proof}

We shall now use the description of 
the collection of sets of semistable points
in terms of GIT-bags in order to study 
basic properties of the corresponding 
system of quotients.
The first statement is the following 
characterization 
of saturated inclusion by means of 
GIT-bags.

\begin{theorem}
\label{saturated}
Let $(D_i,w_i) \in C_T^{\sharp}(X) \cap(K \oplus M)$ 
represent 
two $T$-linearized Weil divisors on $X$.
Then $X^{ss}(D_1,w_1)$ is a $T$-saturated 
subset of $X^{ss}(D_2,w_2)$
if and only if
$\mu(D_1,w_1)^{\circ} \supset \mu(D_2,w_2)^{\circ}$
holds.
\end{theorem}

\begin{proof}
We begin with a preparatory 
observation, characterizing 
closedness of a given $T$-orbit 
in the set $X^{ss}(D_i,w_i)$.

\medskip

\noindent
{\em Claim. }
Consider points $x \in X^{ss}(D_i,w_i)$ and  
$\rq{x} \in W_i := q_X^{-1}( X^{ss}(D_i,w_i))$
such that 
$q_X(\rq{x}) = x$ holds, and
$H \mal \rq{x}$ is closed in 
$\rq{X}$.
Then $T \mal x$ is closed in $X^{ss}(D_i,w_i)$,
if and only if $(H \times T) \mal \rq{x}$ is 
closed in $W_i^{}$. 

\medskip

Let us verify the claim. The ``if'' part is clear.
So, let $T \mal x$ be closed in $X^{ss}(D_i,w_i)$.
Assume that the complement 
$$
Y
\; := \; 
\b{(H \times T) \mal \rq{x}}
\setminus 
(H \times T) \mal \rq{x}
\; \subset \; 
W_i
$$ 
is nonempty. Then Proposition~\ref{goodquotprop}~(ii)
tells us that $x \not \in q_X(Y)$ holds.
On the other hand, we have 
$$
q_X(Y)
\; \subset \;
q_X\bigl(\,\b{(H \times T) \mal \rq{x}}\,\bigr)
\; \subset \; 
\b{q_X(H \times T) \mal \rq{x}}
\; = \;
\b{T \mal x}
\; \subset \; 
X^{ss}(D_i,w_i).
$$
Since $q_X(Y)$ is $T$-invariant, and we assumed $T \mal x$ 
to be closed, this is a contradiction.
Thus, the claim is verified.

We come to the proof of the theorem.
%Set $\sigma_i := \sigma_{H \times T}(D_i,w_i)$.
%Recall that we have 
%$(D_i,w_i) \in \sigma_i^{\circ} \subset \mu(D_i,w_i)^{\circ}$.
First, suppose that $X^{ss}(D_1,w_1)$ is 
$T$-saturated
in $X^{ss}(D_2,w_2)$.
We then have to show
\begin{eqnarray*}
\bigcap_{\omega_1 \in C_T(X); \; (D_1,w_1) \in \omega_1^{\circ}}
\omega_1^{\circ}
& \quad\supset\quad &
\bigcap_{\omega_2 \in C_T(X); \; (D_2,w_2) \in \omega_2^{\circ}}
\omega_2^{\circ}.
\end{eqnarray*}

Consider 
$\omega_1 = \omega_{H \times T}(\rq{x}) \in C_T(X)$ 
as on the left hand side.
Then, by Theorem~\ref{semsitablecrit} and 
Corollary~\ref{clos2face}, 
the orbit $(H \times T) \mal \rq{x}$
is closed in $W_1$,
and the orbit $H \mal \rq{x}$ is closed in
$\rq{X}$.
Thus, $T \mal x$ is closed in $X^{ss}(D_1,w_1)$.
By $T$-saturatedness,
$T \mal x$ is closed in $X^{ss}(D_2,w_2)$.
The above claim tells us that 
$(H \times T) \mal \rq{x}$
is closed in $W_2$.
Thus, Theorem~\ref{semsitablecrit} and 
Corollary~\ref{clos2face}
show $(D_2,w_2) \in \omega_1^{\circ}$.

Conversely, suppose that 
$\mu(D_1,w_1)^{\circ} \supset \mu(D_2,w_2)^{\circ}$
holds.
This clearly implies  
$\mu(D_1,w_1) \ge \mu(D_2,w_2)$.
Theorem~\ref{GIT-bags} thus gives  
$X^{ss}(D_1,w_1) \subset X^{ss}(D_2,w_2)$.
We are left with showing that this is 
a $T$-saturated inclusion.
For this, it suffices to show
that every closed $T$-orbit in 
$X^{ss}(D_1,w_1)$ is as well
closed in $X^{ss}(D_2,w_2)$;
use, for example,
Proposition~\ref{goodquotprop}
and Corollary~\ref{clos2face}.

Consider a closed orbit 
$T \mal x \subset X^{ss}(D_1,w_1)$.
Choose $\rq{x} \in q_X^{-1}(X)$
such that $H \mal \rq{x}$ is closed
in $\rq{X}$. 
By the above claim, 
$(H \times T) \mal \rq{x}$
is closed in $W_1$.
By Theorem~\ref{semsitablecrit} and 
Corollary~\ref{clos2face}, 
the orbit cone
$\omega := \omega_{H \times T}(\rq{x})$
satisfies 
$\omega \in C_T(X)$
and $(D_1,w_1) \in \omega_1^{\circ}$.
The assumption then gives
$(D_2,w_2) \in \omega_1^{\circ}$. 
Using once more Theorem~\ref{semsitablecrit} 
and Corollary~\ref{clos2face}, 
we see that 
$(H \times T) \mal \rq{x}$
is closed in $W_2$.
Consequently, 
$T \mal x$ is closed in $X^{ss}(D_2,w_2)$.
\end{proof}

A a consequence, we can describe the 
{\em qp-maximal\/} $T$-sets of $X$;
these are by definition open 
$T$-invariant subsets $U \subset X$
that admit a good quotient 
$U \to U \quot T$ such that $U \quot T$ 
is quasiprojective and $U$ does not 
occur as a $T$-saturated subset of 
some properly larger $U' \subset X$
admitting a good quotient 
$U' \to U' \quot T$ with $U' \quot T$ 
quasiprojective.

\begin{corollary}
\label{qp-maximal}
Let $\Lambda_T^0(X) \subset \Lambda_T(X)$
consist of all GIT-bags
$\mu_0 \in \Lambda_T(X)$ such that
$\mu_0^{\circ}$ is set theoretically
minimal in $\{\mu^{\circ}; \; \mu \in \Lambda_T(X)\}$.
Then the sets of semistable points
associated to the $\mu_0 \in \Lambda_T^0(X)$
are precisely the qp-maximal $T$-sets
of $X$.
\end{corollary}

\begin{proof}
By Proposition~\ref{qpquotients} 
every qp-maximal $T$-set is the set
of semistable points of a $T$-linearized
Weil divisor on $X$.
Thus, the assertion follows 
from Theorem~\ref{saturated}.
\end{proof}

As the examples discussed in 
the last section of the paper show,
there may exist qp-maximal open subsets
having a non-complete quotient, though
$X$ is assumed to be complete.
For the subcollection of GIT-bags
defining projective quotient spaces,
we obtain a quite simple picture. 

\begin{proposition}
\label{projectivequotients}
Consider the following subcollection of 
the collection $\Lambda_T(X)$ of all 
GIT-bags: 
$$
\Lambda_{T}^\pr(X) 
\; := \; 
\{\mu \in \Lambda_T(X); \;
\forall \ x \in \b{X}\colon \;
\omega_{H \times T}(x)^{\circ} \cap \mu^{\circ} 
\ne \emptyset
 \; \Rightarrow  \;
\omega_{H \times T}(x) \in C_T(X)
\}.
$$  
\begin{enumerate}
\item 
A GIT-bag $\mu \in \Lambda_{T}(X)$
belongs to $\Lambda_{T}^\pr(X)$ if and only 
if the corresponding set of semistable points
has a projective quotient space.
\item 
For any $\mu \in \Lambda^\pr_{T}(X)$,
we have $\mu \in \Sigma_{H \times T}(\b{X})$,
and for any two 
$\mu_1, \mu_2 \in \Lambda_{T}^\pr(X)$,
we have 
$
\mu_1 \le \mu_2
\Leftrightarrow
\mu_1 \preceq \mu_2
$.
\item 
For any two GIT-bags $\mu_1 \in \Lambda_T(X)$
and $\mu_2 \in \Lambda_T^\pr(X)$, we have
$\mu_1 \le \mu_2 \Rightarrow \mu_1 \in \Lambda_T^\pr(X)$.
\end{enumerate}
\end{proposition}

\begin{proof}
Consider a GIT-bag $\mu = \mu(D,w)$ 
and the GIT-cone 
$\sigma := \sigma_{H \times T}(D,w)$.
If we have 
$\mu \in \Lambda_T^\pr(X)$,
then the definition of GIT-bags
and Lemma~\ref{intchamber} yield
$\mu = \sigma$, which is the 
first part of assertion~(ii).
Moreover, Theorem~\ref{semsitablecrit}
yields
$$ 
q_X^{-1}(X^{ss}(D,w))
\; = \; 
\b{X}^{ss}(\sigma).
$$
Since $X^{ss}(D,w) \quot T$ equals
$q_X^{-1}(X^{ss}(D,w)) \quot (H \times T)$,
we infer from the above equation 
that $X^{ss}(D,w) \quot T$ is projective.
Thus, the ``only if'' part of assertion~(i)
is verified.
 
To see the ``if'' part of~(i), suppose that 
the quotient space $X^{ss}(D,w) \quot T$ 
is projective.
Then 
$q_X^{-1}(X^{ss}(D,w)) \quot (H \times T)$
is so. 
Thus, by Lemma~\ref{saturated2},
the inverse image  
$q_X^{-1}(X^{ss}(D,w))$
equals 
$\b{X}^{ss}(\sigma)$.
Consequently, 
Theorem~\ref{semsitablecrit}
gives
\begin{eqnarray*}
\omega_{H \times T}(x) \in C_T(\sigma)
& \Rightarrow &
\omega_{H \times T}(x) \in C_T(X)
\end{eqnarray*}
for all $x \in \b{X}$.
By the definition of a GIT-bag,
this shows $\mu = \sigma$.
Applying once more the above implication,
we obtain 
$\mu \in \Lambda_T^{\pr}(X)$.
So, the proof of assertion~(i)
is complete.

To conclude the proof of~(ii), 
we have to relate two GIT-bags
$\mu_1, \mu_2 \in \Lambda_T^\pr(X)$.
If $\mu_1 \le \mu_2$ holds, then
we obviously have 
$\mu_1 \subset \mu_2$.
Since
$\mu_1$ and $\mu_2$  
are cones of the fan 
$\Sigma_{H \times T}(\b{X})$, 
we have
$\mu_1 \preceq \mu_2$.
Conversely,
$\mu_1 \preceq \mu_2$ 
implies
$\b{X}^{ss}(\mu_1) \supset \b{X}^{ss}(\mu_2)$.
Since we have $\b{X}^{ss}(\mu_i) = q_X^{-1}(X^{ss}(D_i,w_i))$,
where $\mu_i = \mu(D_i,w_i)$, we infer
$\mu_1 \le \mu_2$ from Theorem~\ref{GIT-bags}.

The last assertion is easy to see. 
Let $\mu_i = \mu(D_i,w_i)$.
Then $\mu_1 \le \mu_2$ implies
$X^{ss}(D_1,w_1) \supset X^{ss}(D_2,w_2)$.
Since $X^{ss}(D_2,w_2) \quot T$ is 
projective, and the induced map
$X^{ss}(D_2,w_2) \quot T \to 
X^{ss}(D_1,w_1) \quot T$ is 
dominant, also $X^{ss}(D_1,w_1) \quot T$
must be projective.
This shows $\mu_1 \in \Lambda_T^{\pr}(X)$.  
\end{proof}

Let us indicate, how the description of 
the GIT-equivalence for linearized ample bundle
classes given in~\cite{DoHu} and~\cite{Re}
fits into the present framework. 
For this, recall from Corollary~\ref{amplecone} that
$\kappa_X^{\circ} \subset K_{\QQ}$ is the cone of 
ample divisor classes on $X$.

\begin{proposition}
\label{amplequotients}
Let $(D,w) \in C_T^{\sharp}(X) \cap (K \oplus M)$
represent a $T$-linearized divisor class on $X$,
and consider the corresponding GIT-bag $\mu(D,w)$.  
Then the set of semistable points $X^{ss}(D,w)$ 
arises from an ample $T$-linearized bundle 
if and only if
$\mu(D,w) \in \Lambda_T^0(X)$
and 
$\Pi(\mu(D,w))^{\circ}
\cap \kappa_X^{\circ} 
\ne \emptyset$
hold.
\end{proposition}

\begin{proof}
Suppose first that the set of semistable 
points $X^{ss}(D,w)$ arise from an ample 
$T$-linearized bundle.
Then $X^{ss}(D,w) \quot T$ is
projective, and hence
$X^{ss}(D,w)$ is qp-maximal.
Thus, $\mu(D,w) \in \Lambda_T^0(X)$ holds.
Moreover, according to Theorem~\ref{GIT-bags},
we have $\mu(D,w) = \mu(D',w')$  
with some $(D',w') \in K \oplus M$ 
satisfying 
$D' \in \kappa_X^{\circ}$.
This shows 
$\Pi(\mu(D,w))^{\circ}
\cap \kappa_X^{\circ} 
\ne \emptyset$. 

Conversely, suppose that 
$\mu(D,w) \in \Lambda_T^0(X)$
and 
$\Pi(\mu(D,w))^{\circ}
\cap \kappa_X^{\circ} 
\ne \emptyset$
hold.
Then there is a 
$(D',w') \in \mu(D,w)^{\circ} \cap (K \oplus M)$ 
such that we have $D' \in \kappa_X^{\circ}$. 
The definition of GIT-bags gives
$\mu(D',w')^{\circ} \subset \mu(D,w)^{\circ}$.
Since $\mu(D,w) \in \Lambda_T^0(X)$ holds,
we obtain
$\mu(D',w')^{\circ} = \mu(D,w)^{\circ}$.
This shows $\mu(D',w') = \mu(D,w)$,
and hence Theorem~\ref{GIT-bags} gives the
assertion.
\end{proof}

The preceding two propositions 
allow to rediscover the 
fan structure inside the $T$-ample cone
described in~\cite{DoHu} and~\cite{Re}.
The $T$-ample cone is defined as the cone of 
generated by the $T$-linearized ample 
divisor classes having a nonempty set of 
semistable points, and it is given by
\begin{eqnarray*}
C_T^+(X)
& = &
\Pi^{-1}(\kappa_X)^{\circ}
\cap
C_T^{\sharp}(X).
\end{eqnarray*}

\begin{corollary}
The ample GIT-classes of the $T$-action on
$X$ are in order reversing correspondence 
to the fan of partially open cones 
$\mu \cap C_T^+(X)$, where 
$\mu$ runs through the cones of
$\Lambda_T^0(X)$ 
with
$\kappa_X^\circ \cap \Pi(\mu)^{\circ} \ne
\emptyset$. 
\end{corollary}

\begin{proof}
By Propositions~\ref{amplequotients}
and~\ref{projectivequotients}~(ii)
the partially open cones 
$\mu \cap C_T^+(X)$
as in the assertion are precisely 
the intersections $\sigma \cap C_T^+(X)$,
where $\sigma \in \Sigma_{H \times T}(X)$.
In particular, they form a fan.
The rest is a direct consequence of 
Theorem~\ref{GIT-bags} and 
Proposition~\ref{projectivequotients}~(ii).
\end{proof}

\section{The $\QQ$-factorial case}
\label{sec:Qfactorialcase}

The setup and the notation in this section 
are the same as in the preceding one.
We study the partially ordered collection 
of qp-maximal $T$-sets in terms of GIT-bags 
for the case of a $\QQ$-factorial variety~$X$;
recall that $\QQ$-factoriality means
that $X$ is normal, and that for every
Weil divisor on $X$, some positive
multiple is Cartier.

According to Theorem~\ref{GIT-bags}
and Proposition~\ref{qp-maximal},
the qp-maximal subsets of $X$
are in order reversing
bijection with the GIT-bags in  
$\Lambda^0_T(X) \subset \Lambda_T(X)$,
where $\mu \in \Lambda_T(X)$ belongs 
to $\Lambda_T^0(X)$ if and only if 
$\mu^{\circ}$ is set theoretically 
minimal among the relative interiors 
of all elements of $\Lambda_T(X)$. 
The main result of this section is
the following.

\begin{theorem}
\label{Qfactorial2fan}
Assume that $X$ is $\QQ$-factorial,
and let 
$\mu_1,\mu_2 \in \Lambda_T^0(X)$.
\begin{enumerate}
\item
If $\mu_1 \le \mu_2$ holds, then 
$\mu_1 \preceq \mu_2$ holds, 
and we have 
$$ 
\Star(\mu_1,\mu_2)
\; := \; 
\{\mu \preceq \mu_2; \; \mu_1 \preceq \mu \}
\; \subset \; 
\Lambda_T^0(X).
$$
Moreover, $\mu_1 \le \mu_2$ implies 
$\mu_1 \le \mu \le \mu_2$ for any 
$\mu \in \Star(\mu_1,\mu_2)$.
\item
If there is a 
$\mu_0 \in \Lambda_T^0(X)$
with $\mu_0 \le \mu_1, \mu_2$, 
then we have
$$
\mu_1 \cap \mu_2 \; \preceq \; \mu_1, \mu_2,
\qquad
\mu_1 \cap \mu_2 \; \in \; \Lambda_T^0(X),
\qquad
\mu_0 
\; \le \; 
\mu_1 \cap \mu_2 
\; \le \; 
\mu_1, \mu_2.
$$
\end{enumerate}
\end{theorem}

As a direct application of 
Theorem~\ref{Qfactorial2fan},
we note the following statement on the structure
of the collection of all qp-maximal $T$-sets
of $X$ as a partially ordered set.

\begin{corollary}
\label{minimum}
Assume that $X$ is $\QQ$-factorial.
Then, for any two qp-maximal $T$-sets
$U_1, U_2 \subset X$,
the collection of qp-maximal $T$-sets
$U \subset X$ with $(U_1 \cup U_2) \subset U$ 
is either empty, 
or it contains a unique minimal element.  
\end{corollary}

\begin{proof}
Let $U_1, U_2$ arise from 
GIT-bags
$\mu_1, \mu_2 \in \Lambda_T^0(X)$,
and consider the collection of all
GIT-bags that define sets 
$U \subset X$ of semistable 
points comprising $U_1 \cup U_2$:
$$ 
\Gamma 
\; = \; 
\{\mu_0 \in \Lambda_T^0(X); \; 
\mu_0 \le \mu_1, \; \mu_0 \le \mu_2
\}.
$$
Suppose that $\Gamma \ne \emptyset$ holds.
Then Theorem~\ref{Qfactorial2fan}
gives $\mu_1 \cap \mu_2 \in \Gamma$
and $\mu_0 \le \mu_1 \cap \mu_2$ for any 
$\mu_0 \in \Gamma$.
Hence, the set $X^{ss}(\mu_1 \cap \mu_2)$
is as desired.
\end{proof}

The key to the results in the $\QQ$-factorial 
case is the following observation on the 
ample cone of the variety $X$.

\begin{proposition}
\label{Qample}
The variety $X$ is $\QQ$-factorial if and 
only if the closure
$\kappa_X \subset K_{\QQ}$ of the ample cone
is of full dimension.
\end{proposition}

\begin{proof}
The statement follows from
the well known fact that $\kappa_X$ is of full
dimension in the vector subspace 
$K_{\QQ}^C \subset K_{\QQ}$
generated by the Cartier divisors, 
see~\cite{Kl}. 
\end{proof}

For the proof of Theorem~\ref{Qfactorial2fan}, 
we need a couple of preparatory observations.
The first one holds as well for not necessarily 
$\QQ$-factorial varieties $X$.

\begin{lemma}
\label{notoverlap}
Let $\mu_1,\mu_2 \in \Lambda_T^0(X)$ be
different from each other.
Then we have 
$\mu_1^{\circ} \cap \mu_2^{\circ} = \emptyset$. 
\end{lemma}

\begin{proof}
Suppose that 
$\mu_1^{\circ} \cap \mu_2^{\circ} \ne \emptyset$
holds. 
Then there exists a lattice vector
$(D,w) \in \mu_1^{\circ} \cap \mu_2^{\circ}$.
For the associated GIT-bag, 
we have $\mu(D,w)^{\circ} \subsetneq \mu_i^{\circ}$.
This contradicts $\mu_i \in \Lambda_T^0(X)$.
\end{proof}

\begin{lemma}
\label{claim0}
Assume that $X$ is $\QQ$-factorial,
and consider two orbit cones 
$\omega := \omega_{H \times T}(x)$
and 
$\omega_0 := \omega_{H \times T}(x_0)$
of the $(H \times T)$-action on $\b{X}$
satisfying $\omega_0 \preceq \omega$.
Then 
$\omega_0 \in C_T(X)$
implies  
$\omega \in C_T(X)$.
\end{lemma}

\begin{proof}
Recall that $\omega_0 \in C_T(X)$
merely means 
$\kappa_X^{\circ} \subset \Pi(\omega_0)^{\circ}$.
The statement thus follows from
the fact that $\kappa_X$ is of full dimension.
\end{proof}

\begin{lemma}
\label{claim1}
Assume that $X$ is $\QQ$-factorial.
Let $\sigma \in \Sigma_{H \times T}(\b{X})$ 
and $\omega_0 \in C_T(X)$ such that
$\sigma \cap \omega_0^{\circ} \ne \emptyset$
holds.
Then there is an
$\omega \in C_T(X)$
with 
$\omega_0 \preceq \omega$
and 
$\sigma^{\circ} \subset \omega^{\circ}$. 
\end{lemma}

\begin{proof}
Let $\sigma_0 \preceq \sigma$ be the face
with 
$\sigma_0^{\circ} \cap \omega_0^{\circ} 
\ne \emptyset$.
Then also $\sigma_0$ is a GIT-cone, 
and thus we have 
$\sigma_0 \subset \omega_0$.

By Lemmas~\ref{chamber2semistabset} 
and~\ref{geominterp}, 
any $x_0 \in \b{X}$ with 
$\omega_{H \times T}(x_0) 
= \omega_0$ belongs to $\b{X}^{ss}(\sigma_0)$,
and it has a closed $(H \times T)$-orbit 
inside this set.
So, fix such a point $x_0$, and 
consider the commutative diagram
$$
\xymatrix{
{\b{X}^{ss}(\sigma)}
\ar[r]^{\subset}
\ar[d]
& 
{\b{X}^{ss}(\sigma_0)} 
\ar[d]
\\
{\b{X}^{ss}(\sigma)\quot H \times T }
\ar[r]
& 
{\b{X}^{ss}(\sigma_0)\quot H \times T } 
} 
$$
Since the induced map of quotients is
projective and dominant, it is surjective.
Thus, there is a point
$x \in \b{X}^{ss}(\sigma)$
lying in the same fibre as $x_0$,
and we may even choose $x$ 
such that its $(H \times T)$-orbit
is closed in  $\b{X}^{ss}(\sigma)$.

Consider $\omega := \omega_{H \times T}(x)$.
Lemma~\ref{geominterp} yields
$\sigma^{\circ} \subset \omega^{\circ}$.
Moreover, by Proposition~\ref{goodquotprop}~(ii),
$x_0$ lies in the closure
of the orbit $(H \times T) \mal x$.
Hence, Corollary~\ref{clos2face} gives
$\omega_0 \preceq \omega$.
Lemma~\ref{claim0} then implies
$\omega \in C_T(X)$, and
hence $\omega$ is as desired.
\end{proof}

\begin{proof}[Proof of Theorem~\ref{Qfactorial2fan}]
We first verify the following three 
claims, and then put them together 
to obtain the statements of the theorem.

\medskip

\noindent
{\em Claim 1. }
Let $\mu_1, \mu_2 \in \Lambda_T^0(X)$
such that $\mu_1 \subset \mu_2$ holds.
Then we have $\mu_1 \preceq \mu_2$.

\medskip

To verify this claim, 
let $\nu_1 \preceq \mu_2$ 
denote the face with 
$\mu_1^{\circ} \subset \nu_1^{\circ}$.
Since $\mu_2$ is a union of cones
of the fan $\Sigma_{H \times T}(\b{X})$, 
the cones 
$\tau_{1,k} \in \Sigma_{H \times T}(\b{X})$
with $\tau_{1,k}^{\circ} \subset \nu_1^{\circ}$ 
satisfy
\begin{eqnarray*}
\nu_1^{\circ} 
& = &
\bigcup_k \tau_{1,k}^{\circ}.
\end{eqnarray*} 
 
Since also $\mu_1$ is a union of some of 
the cones $\tau_{1,k}$,
at least one of them, 
say $\tau_{1,0}^{\circ}$,  
satisfies 
$\tau_{1,0}^{\circ}  \subset \mu_1^{\circ}$.
We have to show that 
all $\tau_{1,k}^{\circ}$
are contained in $\mu_1^{\circ}$.
So, suppose that one of them,
say $\tau_{1,1}^{\circ}$, is not. 
Then there must be an orbit cone 
$\omega_0 = \omega_{H \times T}(x_0)$ 
with the following properties:
$$
\omega_0 \in C_T(X),
\qquad
\tau_{1,0}^{\circ} \subset \omega_0^{\circ},
\qquad
\tau_{1,1}^{\circ} \cap \omega_0^{\circ}
= \emptyset.
$$ 

Let  
$\sigma_{1,0} \in \Sigma_{H \times T}(\b{X})$
be a cone
such that $\tau_{1,0} \preceq \sigma_{1,0}$
and 
$\sigma_{1,0}^{\circ} \subset \mu_2^{\circ}$
hold.
Then, according to Lemma~\ref{claim1}, 
there exists an orbit cone 
$\omega = \omega_{H \times T}(x)$
with the following properties:
$$ 
\omega \in C_T(X),
\qquad
\omega_0 \preceq \omega,
\qquad
\sigma_{1,0}^{\circ} \subset \omega^{\circ}.
$$

The last inclusion implies 
$\mu_2^{\circ} \cap \omega^{\circ} \ne \emptyset$.
Since $\mu_2$ belongs to $\Lambda_T^0(X)$,
we can conclude 
$\mu_2^{\circ} \subset \omega^{\circ}$;
otherwise, $\mu_2 \cap \omega$ would
comprise an element $\mu \in \Lambda_T(X)$
with $\mu^{\circ} \subset \mu_2^{\circ}$,
which would be in contradiction to the 
definition of $\Lambda_T^0(X)$.

Thus, we obtained $\mu_2 \subset \omega$.
Consequently, also the face $\nu_1 \preceq \mu_2$
is contained $\omega$.
Since 
$\nu_1^{\circ} \cap \omega_0^{\circ} 
\ne \emptyset$
holds, we can conclude
$\nu_1^{\circ} \subset \omega_0^{\circ}$.
This implies 
$\tau_{1,1}^{\circ} \subset \omega_0^{\circ}$,
a contradiction. 
So, Claim~1 is verified.

\medskip

\noindent
{\em Claim~2.}
Let $\mu_1, \mu_2 \in \Lambda_T^0(X)$
with $\mu_1 \le \mu_2$, and let 
$\mu \preceq \mu_2$ with 
$\mu_1 \preceq \mu$.
Then we have $\mu \in \Lambda_T^0(X)$.

\medskip

Let us check this claim.
Choose any $(D,w) \in \mu^{\circ}$,
and consider the associated 
GIT-bag $\mu(D,w)$.
We show that $\mu(D,w) = \mu$
and $\mu(D,w) \in \Lambda_T^0(X)$ hold. 

First, 
consider any $\omega_2 \in C_T(X)$ with 
$\mu_2^{\circ} \subset \omega_2^{\circ}$.
Let
$\omega_{1} \preceq \omega \preceq \omega_2$ 
denote the faces with 
$\mu_1^{\circ} \subset \omega_1^{\circ}$
and
$\mu^{\circ} \subset \omega^{\circ}$.
Then $\mu_1 \le \mu_2$ implies 
$\omega_1 \in C_T(X)$.
By Lemma~\ref{claim0}, this gives
$\omega \in  C_T(X)$.
Thus, since $(D,w) \in \omega^{\circ}$
holds, 
we obtain $\mu(D,w) \subset \omega$,
and hence $\mu(D,w) \subset \omega_2$.

This consideration shows
$\mu(D,w) \subset \mu_2$.
Since we have 
$(D,w) \in \mu(D,w)^{\circ} \cap \mu^{\circ}$
and $\mu \preceq \mu_2$, we even
obtain $\mu(D,w) \subset \mu$.

To proceed, note that 
there exists a GIT-bag 
$\nu \in \Lambda_T^0(X)$
such that 
$\nu^{\circ} \subset \mu(D,w)^{\circ}$
holds. 
By Claim~1, 
the inclusion 
$\nu \subset \mu_2$ implies
$\nu \preceq \mu_2$. 
Thus, 
$\nu^{\circ} \subset \mu(D,w)^{\circ} \subset \mu^{\circ}$
implies $\nu = \mu$.
Consequently, we obtained
$\mu \in \Lambda_T^0(X)$.

\medskip

\noindent
{\em Claim~3.}
For any three $\mu_1, \mu_2, \mu_3 \in \Lambda_T^0(X)$
with $\mu_1 \preceq \mu_2 \preceq \mu_3$ and 
$\mu_1 \le \mu_3$, 
we have $\mu_1 \le \mu_2 \le \mu_3$.

\medskip

In order to verify 
$\mu_1 \le \mu_2$,
consider $\omega_{2} \in C_T(X)$
with 
$\mu_2^{\circ} \subset \omega_{2}^{\circ}$.
By Lemma~\ref{claim1}, there is an
$\omega_3 \in C_T(X)$ with 
$\omega_{2} \preceq \omega_{3}$
and $\mu_3^{\circ} \subset \omega_3^{\circ}$.
Since we have $\mu_1 \le \mu_3$,
the face $\omega_{1} \preceq \omega_{3}$
with $\mu_1^{\circ} \subset \omega_{1}^{\circ}$
belongs to $C_T(X)$.
Moreover, 
$\mu_1 \preceq \mu_2 \subset \omega_2$ 
and 
$\mu_1^{\circ} \subset \omega_1^{\circ}$
imply $\omega_{1} \preceq \omega_{2}$.
This shows $\mu_1 \le \mu_2$.

Similarly, to see $\mu_2 \le \mu_3$,
consider $\omega_3 \in C_T(X)$ with 
$\mu_3^{\circ} \subset \omega_3^{\circ}$.
For $i=1,2$, let $\omega_i \preceq \omega_3$ 
be the faces with 
$\mu_i^{\circ} \subset \omega^{\circ}$.
Note that $\omega_1 \preceq \omega_2$ 
holds.
Now, $\mu_1 \le \mu_3$ implies
$\omega_1 \in C_T(X)$.
By Lemma~\ref{claim0},
this gives $\omega_2 \in C_T(X)$.
Hence, we can conclude 
$\mu_2 \le \mu_3$,
and Claim~3 is proved.

\medskip

We come to the assertions of the theorem.
For the first one, 
note that $\mu_1 \le \mu_2$ 
implies $\mu_1 \subset \mu_2$, Thus,
Claims~1, 2 and~3 give the desired 
statements.
In the second assertion, 
the case $\mu_1 = \mu_2$ is trivial,
and hence we may assume that 
$\mu_1 \ne \mu_2$ holds.  
Then Lemma~\ref{notoverlap} gives 
$\mu_1^{\circ} \cap \mu_2^{\circ} =
\emptyset$.
Moreover, since $\mu_0 \le \mu_i$ 
implies
$\mu_0 \subset \mu_i$, and hence,
Claim~1 tells us that 
$\mu_0 \preceq \mu_i$ holds.

Consider the faces $\nu_1 \preceq \mu_1$ 
and $\nu_2 \preceq \mu_2$ 
with $(\mu_1 \cap \mu_2)^{\circ} \subset \nu_i^{\circ}$.
Then we have $\mu_0 \preceq \nu_1, \nu_2$.
Thus, Claim~2 yields $\nu_i \in \Lambda_T^0(X)$.
Since we have $\nu_1^{\circ} \cap \nu_2^{\circ}
\ne \emptyset$, Lemma~\ref{notoverlap}
yields $\nu_1 = \nu_2$.
This in turn implies $\mu_1 \cap \mu_2 = \nu_1$.
Thus, $\mu_1 \cap \mu_2$ is a face of 
$\mu_1$ and of $\mu_2$, and we have 
$\mu_1 \cap \mu_2 \in \Lambda_T^0(X)$.
Claim~3 eventually
shows $\mu_0 \le \mu_1 \cap \mu_2 \le \mu_i$.
\end{proof}

We conclude this section with a characterization
of the geometric GIT-quotients in terms of their
describing GIT-bags in the case of a $\QQ$-factorial
variety $X$.
We obtain it as a consequence from the following 
more general statement.

\begin{proposition}
\label{geomquots}
Let $(D,w) \in C_T^{\sharp}(X) \cap (K \oplus M)$,
and consider the associated  
GIT-bag $\mu(D,w)$.  
Then the following 
statements are equivalent:
\begin{enumerate}
\item 
The morphism
$X^{ss}(D,w) \to X^{ss}(D,w) \quot T$
is a geometric quotient.
\item
Any $\omega \in C_T(X)$ with 
$(D,w) \in \omega^{\circ}$ 
satisfies
$ 
\dim(\omega)
=
\dim(\Pi(\omega))
+
\dim(M)
$.
\end{enumerate}
\end{proposition}

\begin{proof}
The quotient $X^{ss}(D,w) \to X^{ss}(D,w) \quot T$
is geometric if and only if all $T$-orbits inside
$X^{ss}(D,w)$ are of full dimension.
The latter holds if and only if 
for all points 
$x \in q_X^{-1}(X^{ss}(D,w))$
with a closed 
$(H \times T)$-orbit
the quotient of isotropy groups
$(H \times T)_x / H_x$
is finite. 
In terms of orbit cones,
this means 
$$ 
\dim(\omega_{H \times T}(x))
\; = \; 
\dim(\omega_H(x)) + \dim(T).
$$
According to Theorem~\ref{semsitablecrit},
the points with a closed $(H \times T)$-orbit
in $q_X^{-1}(X^{ss}(D,w))$ are precisely
those with an orbit cone $\omega \in C_T(X)$ 
satisfying $(D,w) \in \omega^{\circ}$.
This gives the assertion.
\end{proof}

As an immediate consequence of 
Proposition~\ref{Qample},
this characterization of 
geometric quotients breaks 
down in the $\QQ$-factorial case
to the following. 

\begin{corollary}
Assume that $X$ is $\QQ$-factorial,
and let $(D,w) \in C_T^{\sharp}(X) \cap (K \oplus M)$.
Then the quotient 
$X^{ss}(D,w) \to X^{ss}(D,w) \quot T$
is geometric if and only if the   
GIT-bag $\mu(D,w)$ is of full dimension.
\end{corollary}

\section{Examples}
\label{sec:examples}

In this section, we present a couple of examples.
Firstly, we discuss a quite simple example,
a $\KK^*$-action on a Hirzebruch surface,
showing that the intersection of two GIT-bags
need not be a GIT-bag.
Secondly, we treat an ``exotic orbit space''
found by A.~Bia{\l}ynicki--Birula and 
J.~\'{S}wi\c{e}cicka in~\cite[Example~3]{BBSw2};
this is a projective geometric quotient that does 
not arise from an ample bundle.
Finally, we present a non-complete
qp-maximal quotient of a smooth projective
variety.

All our examples are subtorus actions on toric 
varieties $X$. As this setup might be of 
interest for further examples, we briefly explain the
general procedure to obtain the necessary data 
for the study of the GIT-equivalence.  
A toric variety $X$ arises from a fan $\Delta$ 
in the lattice $N_X$ of one parameter subgroups
of the big torus $T_X \subset X$, see~\cite{Fu}.  
As before, we suppose that $X$ is projective, 
and that its divisor class group $\Cl(X)$ is free.

The group $\Cl(X)$ is generated by the classes 
of the invariant prime divisors $D_1, \ldots, D_r$ 
on $X$, which in turn correspond to the 
rays, i.e, the one-dimensional cones, 
$\varrho_1, \ldots, \varrho_r$ of $\Delta$.
By~\cite{Co}, the total coordinate ring 
$\mathcal{R}(X)$ is a polynomial ring in 
$r$ indeterminates, and thus we have 
$\b{X} = \KK^r$ for the corresponding spectrum.
 
In terms of fans, the subset 
$\rq{X} \subset \b{X}$ is obtained 
as follows: let $v_1, \ldots, v_r \in N_X$ 
denote the primitive lattice vectors
generating the rays of $\Delta$,
set $F := \ZZ^r$ and consider the 
linear map $P \colon F \to N_X$ sending the 
$i$-th canonical base vector $e_i \in F$
to $v_i \in N_X$. 
The fan of $\rq{X}$ then consists of 
faces of the positive orthant 
$\delta \subset F_{\QQ}$:
$$ 
\Delta_{\rq{X}}
\; = \; 
\{\rq{\sigma} \preceq \delta; \; P(\rq{\sigma}) \subset \sigma
\text{ for some } \sigma \in \Delta_X\}.
$$

Moreover, the torus $H = \Spec(\KK[\Cl(X)])$ 
acting on $\b{X}$ is the subtorus of $(\KK^*)^r$
having $L := \ker(P)$ as its lattice
of one parameter subgroups.
The canonical map $q_X \colon \rq{X} \to X$
is the toric morphism corresponding to 
the map $P \colon F \to N_X$ of the fans 
$\Delta_{\rq{X}}$ and $\Delta_X$.
Observe that the map $P \colon F \to N_X$
determines a pair of exact sequences,
which are mutually dual to each other:
$$ 
\xymatrix@R-10pt{
0 \ar[r]
&
L \ar[r]
&
F \ar[r]^{P}
&
N_X \ar[r]
&
0
\\
0
&
K \ar[l]
&
E \ar[l]^{Q}
&
M_X \ar[l]
&
0  \ar[l]
}
$$

Note that the subtorus $H \subset (\KK^*)^r$ 
is as well determined
by its {\em weight map\/} $Q \colon E \to K$. 
The weight cone of the $H$-action is given
by $\Omega_H(\b{X}) = Q(\gamma)$, where 
$\gamma \subset E_{\QQ}$ is the positive orthant.  
The fan $\Sigma_{H}(\b{X})$ is the so-called
Gelfan-Kapranov-Zelevinsky decomposition of 
the cone $Q(\gamma)$, that means
that it is the coarsest common refinement of 
all the images $Q(\gamma_0)$, where 
$\gamma_0 \preceq \gamma$, compare~\cite{OdPa}.

The GIT-cone $\kappa_X \in \Sigma_H(\b{X})$
corresponding to $\rq{X} \subset \b{X}$ can be 
calculated as follows, 
compare~\cite[Theorem~10.2]{BeHa2}:
consider the maximal cones 
$\delta_0 \preceq \delta$ 
of the fan $\Delta_{\rq{X}}$,
and determine the corresponding faces
$\gamma_0 = \delta_0^{\perp} \cap \gamma$
of $\gamma$.
Then $\kappa_X$ is the intersection 
over all the images
$Q(\gamma_0)$.
Recall from Corollary~\ref{amplecone}
that the relative interior of 
$\kappa_X$ is the cone of ample divisors
of $X$.    

Now suppose that $T \subset T_X$ is a subtorus
of the big torus of $X$. Then $T \subset T_X$
corresponds to a sublattice $N_T \subset N_X$.
The lifting of the $T$-action to the 
affine multicone $\b{X}$ 
corresponds to an embedding $N_X \to F$ 
with $N_X \cap L = 0$. 
By fixing a lifting of the $T$-action,
we thus decorated the exact sequence
comprising the map $P \colon F \to N_X$
in the following sense:
$$ 
\xymatrix{
0 \ar[r]
&
L \ar[r]
&
F \ar[r]^{P}
&
N_X \ar[r]
&
0
\\
0 \ar[r]
&
L \ar[r]\ar@{=}[u]
&
L \times N_T \ar[r]_{P} \ar[u]
&
N_T \ar[r]\ar[u]
&
0
}
$$

In order to determine the $(H \times T)$-orbit cones
and the fan $\Sigma_{H \times T}(\b{X})$, we have to dualize
the above commutative diagram. The result is the 
following one:  
$$ 
\xymatrix{
0
&
K \ar[l]\ar@{=}[d]
&
E \ar[l]_{Q} \ar[d]^{\rq{Q}}
&
M_X \ar[l] \ar[d]
&
0  \ar[l] 
\\
0
&
K \ar[l]
&
K \oplus M_T \ar[l]
&
M_T \ar[l] 
&
0  \ar[l] 
}
$$
Then the orbit cones of the $(H \times T)$-action on
$\b{X}$ are
precisely the images $\rq{Q}(\gamma_0)$, where 
$\gamma_0 \preceq \gamma$ with $\gamma \subset E_{\QQ}$
being the positive orthant. 

Moreover, the fan 
$\Sigma_{H \times T}(\b{X})$ is 
%again a Gelfand-Kapranov-Zelenvinsky
%decomposition: it is 
the coarsest common refinement of all 
the orbit cones $\rq{Q}(\gamma_0)$.
Finally, the collections $C_T(X)$ and $C_T(\sigma)$ 
for $\sigma \in \Sigma_{H \times T}(\b{X})$ can 
now be directly computed according to their
definitions, and thus it becomes 
possible to determine the collection of 
GIT-bags.

For the computation steps just outlined, it
is most convenient to use suitable computer
programs. For example, we provide a (free)
Maple-Package {\tt TorDiv} doing all the basic
computations needed, see~\cite{BeHaWi}.
In the following examples, we will
therefore omit the computations, and just 
show their results.

\begin{example}
[A $\KK^*$-action on a Hirzebruch surface]
\label{hirzebruch}
As a toric variety, 
the first Hirzebruch surface
$X$ arises from the 
complete fan $\Delta_X$ in $\ZZ^2$ with 
the four rays
$$ 
\varrho_1 := \QQ_{\ge 0} [1,0],
\quad
\varrho_2 := \QQ_{\ge 0} [0,1],
\quad
\varrho_3 := \QQ_{\ge 0} [-1,1],
\quad
\varrho_4 := \QQ_{\ge 0} [0,-1].
$$
We have $\b{X} = \KK^4$, and the action of 
the torus $H = \Spec(\KK[K])$ is given by 
the weight matrix
\begin{eqnarray*}
Q
& = & 
 \left[
\begin{array}{rrrr}
1 & 0 & 1 & 1 \\
0 & 1 & 0 & 1
\end{array}
\right]
\end{eqnarray*} 
Finally, let $T := \KK^*$ act on $\b{X}$ 
with weights $1$, $0$, $-1$, and $0$.
Then the weight matrix of the 
$(H \times T)$-action on $\b{X}$ 
is given by
\begin{eqnarray*}
\rq{Q}
& = & 
 \left[
\begin{array}{rrrr}
1 & 0 & 1 & 1 \\
0 & 1 & 0 & 1 \\
1 & 0 & -1 & 0
\end{array}
\right]
\end{eqnarray*} 

The fans describing the corresponding 
GIT-systems are easy to determine.
Let $w_i \in \ZZ^2$ denote the $i$-th 
column of $Q$. 
Then, the maximal cones of $\Sigma_{H}(\b{X})$
are 
$$
\kappa_1 := \cone(w_1,w_4),
\qquad
\kappa_2 := \cone(w_2,w_4). 
$$ 
Note that the first cone equals $\kappa_X$.
Similarly, denoting by $\rq{w}_i$ the $i$-th 
column of $\rq{Q}$, 
the maximal cones of $\Sigma_{H \times T}(\b{X})$
are
$$
\sigma_1 := \cone(\rq{w}_2,\rq{w}_3,\rq{w}_4),
\qquad
\sigma_2 := \cone(\rq{w}_1,\rq{w}_2,\rq{w}_4),
\qquad
\sigma_3 := \cone(\rq{w}_1,\rq{w}_3,\rq{w}_4). 
$$

All three cones are GIT-bags,
and they even belong to  
$\Lambda_{T}^0(X)$. 
Note that $\sigma_1$ and $\sigma_2$ 
have a 2-dimensional face in common, 
but there is no element in 
$\Lambda_T^0(X)$ which is smaller 
then $\sigma_1$ and $\sigma_2$.
\end{example}

\begin{example}[A.~Bia{\l}ynicki--Birula and 
J.~\'{S}wi\c{e}cicka]
\label{exoticquotients}
Consider the smooth projective variety $X$ 
obtained from
$\PP_{2} \times \PP_{1}$ by blowing 
up first the line ${[z,0,w] \times [0,1]}$
and then the proper transform 
of the line ${[z,w,0] \times [0,1]}$. 
These blow ups are compatible with 
the action of 
$T := \KK^{*} \times \KK^{*}$ on $\PP_{2} \times \PP_{1}$ 
given by
\begin{eqnarray*}
(t_{1}, t_{2}) 
\; \cdot \; 
([x_{0}, x_{1}, x_{2}], [y_{0}, y_{1}]) 
& := & 
([x_{0},
t_{1} x_{1}, t_{1} x_{2}], [y_{0}, t_{2} y_{1}]),
\end{eqnarray*}
Thus there is an induced $T$-action on $X$ 
making the contraction map equivariant. 
We show that there exist precisely four 
different open sets admitting a geometric 
quotient with a projective orbit space,
but only three of them are sets of semistable
points of ample line bundles.

Let us verify these statements.
As announced, we view $X$ 
as a toric variety. 
It arises from the fan $\Delta_{X}$ in 
$N_{X} := \ZZ^{3}$ with the seven rays
$$ 
\begin{array}{llll}
v_{1}:=(1,0,0),& v_{2}:=(0,1,0),& 
v_{3} := (-1,-1,0),& v_{4} := (0,0,1),\\
v_{5} := (0,0,-1),& v_{6} := (1,0,1),& 
v_{7} :=(0,1,1),
\end{array}
$$
and the ten maximal cones ---
we denote for short by $C_{ijk}$ the cone 
generated by $v_{i},v_{j}, v_{k}$:
$$ 
C_{2 3 5},\; C_{3 4 7},\; C_{2 3 7},\; 
C_{1 3 5},\; C_{3 4 6},\; C_{1 3 6},\; 
C_{1 2 5},\; C_{4 6 7},\; C_{2 6 7},\; C_{1 2 6}.
$$

The spectrum of the total coordinate 
ring $\mathcal{R}(X)$ is $\b{X}=\KK^7$,
and the weight map for the action of 
$H = \Spec(\KK[K])$ on $\b{X}$ 
is given (w.r. to the canonical bases) by
\begin{eqnarray*}
Q 
& = &
\left[ \begin{array}{ccccccc} 1&1&1&0&0&0&0\\
\noalign{\medskip}0&0&0&1&1&0&0\\
\noalign{\medskip}-1&0&0&0&1&1&0\\
\noalign{\medskip}0&-1&0&0&1&0&1
\end {array}\right].
\end{eqnarray*}
According to the preceding remarks, 
the closure $\kappa_X \subset K_{\QQ} = \Cl(X)$
of  the ample cone is given by
\begin{eqnarray*}
\kappa_{X} 
& = & 
\cone(
(1, 0, 0, 0), \;
(2, 1, 0, 0), \; 
(0, 1, 1, 1), \; 
(1, 1, 0,1)).
\end{eqnarray*}

The sublattice $N_{T} \subset N_{X}$ 
describing the $T$-action on $X$ then 
is generated by the vectors 
$(1,1,0)$ and $(0,0,1)$.
Hence, we may work with
\begin{eqnarray*}
\rq{Q} 
& = & 
\left[ \begin{array}{ccccccc}
1&1&1&0&0&0&0\\
\noalign{\medskip}0&0&0&1&1&0&0\\
\noalign{\medskip}-1&0&0&0&1&1&0\\
\noalign{\medskip}0&-1&0&0&1&0&1\\
\noalign{\medskip}1&1&0&0&0&0&0\\
\noalign{\medskip}0&0&0&1&0&0&0
\end{array} \right]. 
\end{eqnarray*}

Now, one has to compute the fan 
$\Sigma_{H\times T}(\b{X})$ 
associated to the $(H \times T)$-action 
on $\b{X}$.
It 
%turns out that $\Sigma_{H\times T}(\b{X})$ 
has the vectors
$$ 
\begin{array}{llll}
w_{1} := (1,0,0,0,0,0), \quad 
& 
w_{2} := (0,1,0,0,0,1),  \quad 
& 
w_{3} := (0,1,1,1,0,0,0), 
\\
w_{4} := (1,0,-1,0,1,0),  \quad 
& 
w_{5} := (1,0,0,-1,1,0),  \quad 
& 
w_{6} := (1,0,0,0,1,0), 
\\
w_{7} := (0,0,1,0,0,0),  \quad 
& 
w_{8} := (0,0,0,1,0,0)
\end{array}
$$
as the primitive generators of its rays,
and it has precisely four full-dimensional 
cones, namely:
$$
\begin{array}{ll}
%\mu_{1} & := 
\cone(w_{1} , w_{2} , w_{3}, w_{4}, w_{5}, w_{6}), \quad
&
%\mu_{2} & := 
\cone(w_{1} , w_{2} , w_{3}, w_{7}, w_{5}, w_{6}), \\
%\mu_{3} & := 
\cone(w_{1} , w_{2} , w_{3}, w_{4}, w_{8}, w_{6}), \quad
&
%\mu_{4} & := 
\cone(w_{1} , w_{2} , w_{3}, w_{7}, w_{8}, w_{6}).
\end{array}
$$

It turns out that these cones are 
precisely the GIT-bags of full dimension.
Proposition~\ref{geomquots} thus tells us
that there are precisely four different
geometric quotients arising 
from linearized bundles.
Moreover, by Proposition~\ref{projectivequotients},
the associated quotient spaces are projective. 
Finally, Proposition~\ref{amplequotients} 
yields that the second one of the above GIT-bags
describes a quotient that 
does not arise from an ample line bundle.  
\end{example}

\begin{example}
\label{maxnotproj}
We present a smooth projective
variety $X$ of dimension three and
a $\KK^*$-action on $X$ that
has qp-maximal sets with non-complete
quotient spaces.

Our $X$ is the toric variety 
arising from the fan $\Delta_X$  
in $N_{X} := \ZZ^{3}$,
which has the vectors
$$ 
v_{1} :=(1,0,0), \; v_{2}:=(0,1,0), \; 
v_{3} := (-1,0,1), \; v_{4} := (0,-1,1), \; 
v_{5} := (0,0,-1)
$$
as the primitive generators of its rays,
and which has the following list of 
maximal cones:
$$ 
\begin{array}{lll}
C_{1} := \cone(v_{1},v_{2},v_{3}),& 
C_{2} := \cone(v_{1},v_{3},v_{4}), & 
C_{3} := \cone(v_{1},v_{2},v_{5}), \\
C_{4} := \cone(v_{2},v_{3},v_{5}), & 
C_{5} := \cone(v_{3},v_{4},v_{5}), &
C_{6} := \cone(v_{1},v_{4},v_{5}).
\end{array}
$$ 

We have $\b{X} = \KK^5$, and the action
of the torus $H = \Spec(\KK[K])$ on $\b{X}$ is given by 
the weight map 
\begin{eqnarray*}
Q 
& := & 
\left[ \begin{array}{ccccc} 
-1&1&-1&1&0 \\ 
%\noalign{\medskip} 
0&1&0&1&1 
\end{array} \right].
\end{eqnarray*}
Moreover, the closure
$\kappa_{X} \subset K_{\QQ} = \Cl_{\QQ}(X)$ 
of the ample cone of $X$ is
generated by the vectors 
$(1,1)$ and $(0,1)$.

Now, consider $\KK^{*}$-action on $X$ corresponding
to the sublattice $N_{T}$ of $N_{X}$ 
generated by $(2,-4,1)$.
A lifting of this action to $\b{X}$ is given by 
the weight matrix
\begin{eqnarray*}
\rq{Q} 
& := & 
\left[ \begin{array}{ccccc} 
-1&1&-1&1&0 \\ 
\noalign{\medskip} 0&1&0&1&1 \\ 
\noalign{\medskip} 3&-4&1&0&0 
\end{array} \right].
\end{eqnarray*}

The fan $\Sigma_{H \times T}(\b{X})$ 
lives in the three-dimensional 
lattice $K \oplus M_{T}$, and
the cone $\sigma$ generated by
$(-1,1,2)$ is an element of 
$\Sigma_{H \times T}(\b{X})$.
It turns out that $\sigma$ is a 
GIT-bag. 
Hence, by Corollary~\ref{qp-maximal},
the associated set of semistable points
$X(\mu)$ is qp-maximal. 

However, the criterion for projectivity
given in Proposition~\ref{projectivequotients}
is not fulfilled: the orbit cone 
$$ 
\omega 
\; := \: 
\cone((-1,0,3), (-1,0,1), (0,1,0))
$$ 
contains $\sigma^{\circ}$ 
in its relative interior, 
but the image of $\omega^{\circ}$ 
in $K$ does not intersect $\kappa_{X}^{\circ}$.
Therefore, $X(\mu) \quot T$ is not
projective.
\end{example}

\end{document}